\documentclass[10pt,a4paper]{article}
\usepackage{hyperref}
\usepackage{amsmath}
\usepackage{amsfonts}
\usepackage{amsthm}
\usepackage{pst-plot}
\newtheorem{theorem}{Theorem}[section]
\newtheorem{lemma}[theorem]{Lemma}
\newtheorem{cor}[theorem]{Corollary}
\newtheorem{proposition}[theorem]{Proposition}
\theoremstyle{remark}
\newtheorem{remark}[theorem]{Remark}
\theoremstyle{definition}
\newtheorem{definition}[theorem]{Definition}
\newtheorem{assumption}[theorem]{Assumption}
\numberwithin{equation}{section}
\hypersetup{pdftitle={Markovianity and ergodicity for a surface growth PDE},
pdfsubject={This paper studies a model in surface growth. We provide the existence
of a weak solution satisfying energy inequalities and the Markov property. Furthermore,
we establish that any such solution is strong Feller and has a unique invariant measure.},
pdfauthor={D. Bloemker, F. Flandoli, M. Romito},
pdfkeywords={Surface growth model, weak energy solutions, Markov solutions, strong Feller property, ergodicity.}}
\begin{document}
\title{Markovianity and ergodicity for a surface growth PDE}
\author{D.~Bl\"omker\thanks{Institut f\"ur Mathematik,  Universit\"at Augsburg, 86135 Augsburg, Germany.
                            E-mail address: \texttt{dirk.bloemker@math.uni-augsburg.de}},
    \and F.~Flandoli\thanks{Dipartimento di Matematica Applicata, Universit\`a di Pisa, via Bonanno Pisano 25/b, 56126 Pisa, Italia.
                            E-mail address: \texttt{flandoli@dma.unipi.it}},
    \and M.~Romito\thanks{Dipartimento di Matematica, Universit\`a di Firenze, Viale Morgagni 67/a, 50134 Firenze, Italia.
                            E-mail address: \texttt{romito@math.unifi.it}}}
\date{}
\maketitle
\begin{abstract}
The paper analyses a model in surface growth, where uniqueness
of weak solutions seems to be out of reach. 
We provide the existence of a weak martingale solution satisfying energy
inequalities and having the Markov property. Furthermore,
under non-degeneracy conditions on the noise, we establish
that any such solution is strong Feller and has a unique
invariant measure.
\smallskip\par
{\footnotesize
\par}
\end{abstract}
\section{Introduction}
The paper deals with a model arising in the theory of growth of surfaces,
where an amorphous material is deposited in high vacuum on an initially
flat surface. Details on this model can be found in Raible et al.\ 
\cite{Ra-Ma-Li-Mo-Ha-Sa:00}, \cite{Ra-Li-Ha:00} or
Siegert \& Plischke \cite{Si-Pl:94:2nd}. After rescaling
the equation reads
\begin{equation}\label{e:dirkeq}
\dot h=-h_{xxxx}-h_{xx}+(h_x^2)_{xx}+\eta
\end{equation}
with periodic boundary conditions on the interval $[0,L]$, where
the noise $\eta$ is white in space and time.

Periodic boundary conditions are the standard condition in
these models. Sometimes the model has been considered also
on the whole real line, even though we do not examine this case.
We remark that from a mathematical point of view Neumann
or Dirichlet boundary conditions are quite similar for the
problem studied here. The key point ensured by any of these
boundary conditions is that there is a suitable cancellation
in the non-linearity, namely
$$
\int_0^L h\,(h_x^2)_{xx}\,dx =0,
$$
which is the main (and only) ingredient to 
derive useful a-priori estimates.

The main terms are the dominant linear operator, and the
quadratic non-linearity. The linear instability $-h_{xx}$,
which leads to the formation of hills, is sometimes neglected
(as we shall do in the analysis of the long time behaviour in
Section \ref{s:ergodic}). 

For general surveys on surface growth processes and molecular beam
epitaxy see Barab\'asi \& Stanley \cite{BaSt95} or Halpin-Healy \& Zhang
\cite{HaZh95}. Recently the equation has also become a model for
ion-sputtering, too, where a surface is eroded by a ion-beam,
see Cuerno \& Barab\'asi \cite{CuBa95}, Castro et al.\ \cite{CCVG05}.

Sometimes one adds to the model an additional non-linear term 
$-h_x^2$ of Kuramoto-Sivashinsky type, but in the present form
the equation is mass conserving (i.e.\ $\int_0^Lh(t,x)\,dx=0$),
as $h$ is the height subject to a moving frame, where the mean
growth of the surface is scaled away.
\subsubsection*{Known results on the model}
Before stating the main results of the paper, we give a short
account of the previously known results concerning both the
deterministic and the stochastic version of the model.
\begin{itemize}
\item[${\scriptstyle\triangleright}$] If $\eta=0$ then the equation
   has an absorbing set in $L^2$, although the solution
   may not be unique. (Stein \& Winkler \cite{StWiP}).
\item[${\scriptstyle\triangleright}$] There exists a unique local
   solution in $L^p([0,\tau),H^1)\cap C^0((0,\tau),H^1)$ for initial
   conditions in $H^s$ with $s>1-\frac1p$ and $p>8$ (see Bl\"omker \& Gugg \cite{BlGu}).
\item[${\scriptstyle\triangleright}$] There are stationary solutions,
   which can be constructed as limit points of stationary solutions
   of Galerkin approximations (see Bl\"omker \& Hairer \cite{BlHa04}).
\item[${\scriptstyle\triangleright}$] There are weak martingale
   solutions by means of Galerkin approximation (see Bl\"omker,
   Gugg, Raible \cite{BlGu02, BlGuRa02}).
\end{itemize}
The main problem of the model, which is shared by both the deterministic
and the stochastic approach, is the lack of uniqueness for weak solutions.
This is very similar to the celebrated Navier-Stokes equation.
With this problem in mind, a possible approach to analyse the
model is to look for solutions with \emph{special} properties,
possibly with a physical meaning, such as the balance of energy
-- we shall often refer to it as \emph{energy inequality} -- or
the Markov property.
\subsubsection*{Main results}
Here we use the method developed by Flandoli \& Romito \cite{FlRo06a},
\cite{FlRo06b} and \cite{FlRo06c} in order to establish the
existence of weak solutions having the Markov property. 
For the precise formulation of the concept of solution 
see Definitions \ref{d:weakms} and \ref{d:energyms}.

The method is essentially based on showing a multi-valued
version of the Markov property for sets of solutions and then
applying a clever selection principle (Theorem \ref{t:exmarkov}).
The original idea is due to Krylov \cite{Kry} (see also
Stroock \& Varadhan \cite[Chapter 12]{StVa}).

A key point in this analysis is the definition of
weak martingale solutions. The above described
procedure needs to handle solutions which incorporate
all the necessary bounds on the size of the process
(solution to the SPDE) in different norms. These
bounds must be compatible with the underlying
Markov structure. This justifies the extensive
study of the energy inequality in Section
\ref{s:marprob}.

Once the existence of at least one Markov family
of solutions is ensured, the analysis of such
solutions goes further. Indeed, the selection
principle provides a family of solutions whose
dependence with respect to the initial conditions
is just measurability. By slightly restricting
the set of initial condition, this dependence
can be improved to continuity in the total
variation norm (or \emph{strong Feller} in
terms of the corresponding transition
semigroup). In few words, we show that the
smaller space $H^1_{per}$  (see next section
for its precise definition) is the natural
framework for the stochastic model.

Our last main result concerns the long time
behaviour of the model. We are able to show
that any Markov solution has a unique invariant
measure whose support covers the whole state
space. In principle the existence of stationary
states has been already proved by Bl\"omker
\& Hairer \cite{BlHa04}. Their result are not
useful in this framework, as we have a transition
semigroup that depends on the generic selection
under analysis and in general\footnote{Unless
the problem is well-posed!} is not obtained by
a suitable limit of Galerkin approximations.
In this way, our results are more powerful,
as they apply to \emph{every} Markov solution.
The price to pay is that the proof of existence
of an invariant measure is painfully long
and technical (see Section \ref{s:ergodic}).
 
We finally remark that, even though we know by these
results that every Markov solution is strong Feller
and converges to its own invariant measure, well
posedness is still an open problem for this model
and these result essentially do not improve our
knowledge on the problem. Even the invariant
measures are different, as they depend from
different Markov semigroups.
\subsubsection*{A comparison with previous results on the Markov property}
There are several mathematical interests in this model, in
comparison with the theory developed in Flandoli \& Romito
\cite{FlRo06a}, \cite{FlRo06b} and \cite{FlRo06c} for the
Navier-Stokes equations.
Essentially, in this model we have been able to find
the \emph{natural space} for the Markov dynamics, thus
showing the existence of the (unique) invariant measure.
This is still open, in the framework of Markov selections,
for the Navier-Stokes equations.

Another challenge of this model has concerned the analysis
of the energy inequality. Here the physics of the model
requires a noise white in time and space, while the
analysis developed in the above cited papers has been
based on a trace-class noise with quite regular
trajectories.

Finally, we remark that there is a different approach to handle
the existence of solutions with the Markov property, based on
spectral Galerkin methods, which has been developed by Da Prato
\& Debussche \cite{DPDe} (see also Debussche \& Odasso \cite{DeOd})
for the Navier-Stokes equations (no result with these techniques
is known on the model analysed in this paper). Their methods are
similar to \cite{BlGuRa02, BlGu02, BlHa04}.
\subsubsection*{Layout of the paper}
The paper is organised as follows. In Section \ref{s:marprob}
we state the martingale problem and define weak and energy
solutions. We also give a few restatements of the energy balance.
We next show in Section \ref{s:exist} that there is at least
one family of energy solutions with the Markov property.
In Section \ref{s:sfeller} we show that the transition semigroup
associated to any such solution has the strong Feller property.
Existence and uniqueness of the invariant measure is then
shown in Section \ref{s:ergodic}.

Finally, Sections \ref{sec:WMS} and \ref{s:tools} contain
a few technical results that are used along the paper.
They have been confined in the last part of the
paper to ease the reader from such details
and focus on the main topics.
\section{The martingale problem}\label{s:marprob}
\subsection{Notations and assumptions}
Let $\mathcal{D}^\infty$ be the space of infinitely differentiable
$L$-periodic functions on $\mathbf{R}$ with zero mean in $[0,L]$.
We work with periodic boundary conditions on $[0,L]$ and mean zero
and we define
$$
L^2_{per}=\{h\in L^2(0,L):\int_0^Lh(x)\,dx=0\}.
$$
the spaces $H^1_{per}$, $H^2_{per}$, etc.\ are defined similarly
(see for example \cite[Section 2]{BlGu02}).

Let $\Delta$ be the operator $\partial_x^2$ on $L^2_{per}$
subject to periodic boundary conditions. The leading linear 
operator in \eqref{e:dirkeq} is  $A=-\Delta^2$.
Let $(e_k)_{k\in\mathbf{N}}$ be the orthonormal basis of $L^2_{per}$ given
by the trigonometric functions and let $\lambda_k$ be the eigenvalues of $A$ such that
$$
Ae_k=\lambda_ke_k.
$$
Notice that $\lambda_k\sim -k^4$.

Let $\mathcal{Q}:L^2_{per}\to L^2_{per}$ be a bounded linear operator such that
$$
\mathcal{Q}e_k=\alpha_k^2e_k,
\qquad k\in\mathbf{N},
$$
so that $\mathcal{Q}$ is non-negative self-adjoint operator.
This is sufficient to model all kinds of spatially homogeneous Gaussian noise $\eta$
such that
$$
\mathbb{E}\eta(t,x)=0
\qquad \mathrm{and}\qquad
\mathbb{E}\eta(t,x)\mathbb{E}\eta(s,y)=\delta(t-s)q(x-y)
$$
where $q$ is the the spatial correlation function (or distribution).
Now $\mathcal{Q}=q\star$, which is the convolution operator with $q$.
For details see Bl\"omker \cite{Bl-noise} and the references therein.

In a formal way we can rewrite \eqref{e:dirkeq} as an abstract 
stochastic evolution equation
$$
dh=(Ah-\Delta h + B(h,h))\,dt +dW,
$$
where $W$ is a suitable $\mathcal{Q}$-Wiener process (for details see \eqref{e:defW}), 
and $B(u,v)= -\Delta(\partial_x u\cdot \partial_x v).$ 
\subsubsection{The underlying probability structure}
Let $\Omega=C([0,\infty);H^{-4}_{per})$ and let $\mathcal{B}$
be the $\sigma$-algebra of Borel subsets of $\Omega$.
Let $\xi:\Omega\to H^{-4}_{per}$ be the \emph{canonical
process} on $\Omega$, defined as $\xi(t,\omega)=\omega(t)$.

For each $t\ge0$, let $\mathcal{B}_t=\sigma[\xi(s):\ 0\le s\le t]$
be the $\sigma$-field of events up to time $t$ and
$\mathcal{B}^t=\sigma[\xi(s):\ s\ge t]$ be the $\sigma$-field
of events after time $t$. The $\sigma$-field $\mathcal{B}_t$
can be seen as the Borel $\sigma$-field of
$\Omega_t=C([0,t];H^{-4}_{per})$ and, similarly, $\mathcal{B}^t$
as the Borel $\sigma$-field of $\Omega_t=C([t,\infty];H^{-4}_{per})$.
Notice that both $\Omega_t$ and $\Omega^t$ can be seen as Borel
subsets of $\Omega$ (by restriction to corresponding sub-intervals).
Define finally the \emph{forward shift}\, $\Phi_t:\Omega\to\Omega^t$,
defined as
\begin{equation}\label{e:Phidef}
\Phi_t(\omega)(s)=\omega(s-t),
\qquad s\ge t.
\end{equation}

Given a probability measure $P$ on $(\Omega,\mathcal{B})$ and $t>0$, we
shall denote by $\omega\mapsto P|_{\mathcal{B}_t}^\omega:\Omega\to\Omega^t$
a regular conditional probability distribution of $P$ given $\mathcal{B}_t$%
\footnote{Notice that $\Omega$ is a Polish space and $\mathcal{B}_t$ is
countably generated, so a regular conditional probability distribution
does exist and is unique, up to $P$-null sets.}. In particular,
$P|_{\mathcal{B}_t}^\omega[\omega':\ \xi(t,\omega')=\omega(t)]=1$ and,
if $A\in\mathcal{B}_t$ and $B\in\mathcal{B}^t$, then
$$
P[A\cap B]=\int_A P|_{\mathcal{B}_t}^\omega[B]\,P(d\omega).
$$
One can see the probability measures $(P|_{\mathcal{B}_t}^\omega)_{\omega\in\Omega}$ as
measures on $\Omega$ such that
$P|_{\mathcal{B}_t}^\omega[\omega'\in\Omega:\omega'(s)=\omega(s),\text{ for all }s\in[0,t]]=1$
for all $\omega$ in a $\mathcal{B}_t$-measurable $P$-full set.
We finally define the \emph{reconstruction} of probability
measures (details on this can be found in Stroock \& Varadhan
\cite[Chapter 6]{StVa}).
\begin{definition}\label{e:reconstr}
Given a probability measure $P$ on $(\Omega,\mathcal{B})$, $t>0$
and a $\mathcal{B}_t$-measurable map $Q:\Omega\to\Pr(\Omega^t)$
such that $Q_\omega[\xi_t=\omega(t)]=1$ for all $\omega\in\Omega$,
$P\otimes_t Q$ is the unique probability measure on $(\Omega,\mathcal{B})$
such that
\begin{enumerate}
\item[1.] $P\otimes_t Q$ agrees with $P$ on $\mathcal{B}_t$,
\item[2.] $(Q_\omega)_{\omega\in\Omega}$ is a regular conditional
probability distribution of $P\otimes_t Q$, given $\mathcal{B}_t$.
\end{enumerate}
\end{definition}
\subsection{Solutions to the martingale problem}
\begin{definition}[weak martingale solution]\label{d:weakms}
Given $\mu_0\in\Pr(L^2_{per})$, a probability measure $P$ on
$(\Omega,\mathcal{B})$ is a solution, starting at $\mu_0$,
to the martingale problem associated to equation \eqref{e:dirkeq}
if
\begin{enumerate}
\item[\rm\footnotesize\bfseries\textsf{[W1]}]
$P[L^2_{loc}([0,\infty);H^1_{per})]=1$,
\item[\rm\footnotesize\bfseries\textsf{[W2]}] for every
$\varphi\in\mathcal{D}^\infty$, the process $(M_t^\varphi,\mathcal{B}_t,P)_{t\ge0}$,
defined $P$-a.\ s.\ on $(\Omega,\mathcal{B})$ as
$$
M_t^\varphi=\langle\xi(t)-\xi(0),\varphi\rangle
           +\int_0^t\langle\xi(s),\varphi_{xxxx}+\varphi_{xx}\rangle\,ds
           -\int_0^t\langle(\xi_x(s))^2,\varphi_{xx}\rangle\,ds
$$
is a Brownian motion with variance $t|\mathcal{Q}^{\frac12}\varphi|_{L^2}^2$,
\item[\rm\footnotesize\bfseries\textsf{[W3]}] the marginal at time $0$ of $P$ is $\mu_0$
\end{enumerate}
\end{definition}
\begin{remark}\label{r:equiv}
It is not difficult to prove that the definition of weak martingale
solution given above coincides with the usual definition given
in terms of existence of an underlying probability space and
a Wiener process. This same equivalence is proved in Flandoli
\cite{FlaCIME} for the Navier-Stokes equations and one can
proceed similarly in this case.
\end{remark}
Define, for every $k\in\mathbf{N}$, the process
$\beta_k(t)=\frac1{\alpha_k}M^{e_k}_t$ (and
$\beta_k=0$ if $\alpha_k=0$). Under any weak
martingale solution $P$, the
$(\beta_k)_{k\in\mathbf{N}}$ are a sequence
of independent one-dimensional standard Brownian motions.

Similarly, the process
\begin{equation}\label{e:defW}
W(t)=\sum\alpha_k\beta_k(t)e_k
\end{equation}
is, under any weak martingale solution $P$, a $\mathcal{Q}$-Wiener process and the
process\footnote{The process $Z$ can be equivalently defined as
$$
Z(t,\omega)
=W(t,\omega)+\int_0^tA\mathrm{e}^{A(t-s)}W(s,\omega)\,ds.
$$
The process $Z$ is defined, in some sense, \emph{path-wise} and so versions
of this process cannot be used.}
\begin{equation}\label{e:defZ}
Z(t)
=\sum_{k\in\mathbf{N}}\alpha_k\int_0^t\mathrm{e}^{(t-s)\lambda_k}\,d\beta_k(s)e_k
\end{equation}
is the associated Ornstein-Uhlenbeck process starting at $0$. 
The sum above has to be understood as the limit in $L^2(\Omega)$, and we
know that, under any weak martingale solution, it converges.

Notice that, obviously, $Z$ and $W$ are random variables
on $\Omega$. We state a first regularity result for $Z$,
that we shall use in the definition below.
\begin{lemma}\label{l:Zreg}
Given a weak martingale solution $P$, then for every $T>0$, $p\ge1$, and $s\in[0,3/2)$
$$
Z\in L^p(\Omega\times(0,T);W^{1,4}_{per}),
$$
as for some $\lambda>0$,
$$
\sup_{T>0}\frac1T \mathbb{E}\bigl[\int_0^T \exp\{\lambda |Z(t)|^2_{W^{1,4}}\}\,dt\bigr]<\infty.
$$
Furthermore,
$$
Z\in L^p(\Omega, L^\infty([0,T],H^s_{per})).
$$
Due to $Z\in\Omega$, we thus have $Z$ is a.\ s.\ weakly continuous with values in $H^s_{per}$.
\end{lemma}
This lemma will be proved in Lemmas \ref{lem:ZregW14} and \ref{lem:ZregLinfty} in Section \ref{sec:WMS}.
\begin{definition}[energy martingale solution]\label{d:energyms}
Given $\mu_0\in\Pr(H)$, a probability measure $P$ on $(\Omega,\mathcal{B})$
is an \emph{energy martingale solution} to equation \eqref{e:dirkeq}
starting at $\mu_0$ if
\begin{enumerate}
\item[\rm\footnotesize\bfseries\textsf{[E1]}] $P$ is a weak martingale solution starting at $\mu_0$,
\item[\rm\footnotesize\bfseries\textsf{[E2]}] $P[V\in L^\infty_{loc}([0,\infty);L^2_{per})\cap L^2_{loc}([0,\infty);H^2_{per})]=1$,
\item[\rm\footnotesize\bfseries\textsf{[E3]}] there is a set   $T_P\subset(0,\infty)$ of null Lebesgue
  measure such that for all $s\not\in T_P$ and all $t\ge s$,
$$
P[\mathcal{E}_t(V,Z)\le\mathcal{E}_s(V,Z)]=1,
$$
\end{enumerate}
where $V(t,\omega)=\xi(t,\omega)-Z(t,\omega)$, for $t\ge 0$, and the energy functional
$\mathcal{E}$ is defined as
$$
\mathcal{E}_t(v,z)
= \frac12|v(t)|_{L^2}^2
 +\int_0^t(|v_{xx}|_{L^2}^2-|v_x|_{L^2}^2-\langle v_x,z_x\rangle_{L^2}-\langle 2v_x z_x+(z_x)^2,v_{xx}\rangle_{L^2}).
$$
\end{definition}
\begin{remark}[The equation for $V$]\label{r:eqweak}
Let $P$ be an \emph{energy martingale solution}, then it is easy
to see that, by definition, $M^\varphi_t=\langle W(t),\varphi\rangle$
for all $\varphi\in\mathcal{D}^\infty$. Moreover,
$$
\langle Z(t),\varphi\rangle+\int_0^t\langle Z(s),\varphi_{xxxx}\rangle\,ds
=\langle W(t),\varphi\rangle,
$$
and thus
$$
\langle V(t)-\xi(0),\varphi\rangle
+\int_0^t(\langle V,\varphi_{xxxx}+\varphi_{xx}\rangle
         +\langle Z - (V_x+Z_x)^2,\varphi_{xx}\rangle)\,ds=0,
$$
or, in other words, $V$ is a weak solution (i.e.\ in the sense of distributions) 
to the equation,
$$
\dot V+V_{xxxx}+V_{xx}+Z_{xx}=[(V_x+Z_x)^2]_{xx},
$$
with initial condition $V(0)=\xi(0)$.
\end{remark}
\begin{remark}[Finiteness of the energy]
Given a \emph{energy martingale solution} $P$, we aim to
show that, under $P$, the energy $\mathcal{E}_t$ is almost
surely finite. Indeed, by {\rm\textbf{\footnotesize\textsf{[E2]}}},
it follows that $V(t)$ is $P$-a.\ s.\ weakly continuous
in $L^2_{per}$ (see for example Lemma 3.1.4 of Temam \cite{Tem}),
and so the function $|V(t)|_{L^2}^2$ is defined point-wise
in the energy estimate. Similarly, the other terms are
also $P$-a.\ s.\ finite by {\rm\textbf{\footnotesize\textsf{[E2]}}}
and the regularity properties of $Z$ under $P$
(Lemma \ref{l:Zreg}).
\end{remark}
\begin{remark}[Measurability of the energy and equivalent formulations]
This last remark is concerned with the measurability
issues related to the energy inequality and with some
equivalent formulations of property
{\rm\textbf{\footnotesize\textsf{[E3]}}} of the above
definition. We first prove in the next lemma that property
{\rm\textbf{\footnotesize\textsf{[E3]}}} is quite strong
and that, in a sense that will be clarified below, the
energy inequality is an intrinsic property of the solution
to the original problem \eqref{e:dirkeq}, and does not depend
on the splitting $V+Z$. A similar result was proved in Romito
\cite{Rom01} for the Navier-Stokes equations. We then show
measurability of the energy balance functional and give some
equivalent formulations of the energy inequality.
\end{remark}
Before stating the lemma, we introduce some notations. Let
$z_0\in H^1_{per}$ and $\alpha\ge0$, and let $\widetilde{Z}=\widetilde{Z}_{\alpha,z_0}$ be the solution
to
\begin{equation}\label{e:ztilde}
\dot{\widetilde{Z}}=-{\widetilde{Z}}_{xxxx}-\alpha\widetilde{Z}+\eta,
\qquad\widetilde{Z}(0)=z_0.
\end{equation}
The process $\widetilde{Z}$ is given by $\widetilde{Z}=Z+w$, where $w$ solves the
(deterministic) problem
\begin{equation}\label{e:ztildemenoZ}
\dot w=-w_{xxxx}-\alpha\widetilde{Z},
\qquad w(0)=z_0,
\end{equation}
and so it is well defined $P$-a.\ s., for every martingale solution
$P$. Define suitably $\widetilde{V}=\widetilde{V}_{\alpha,z_0}$ as
${\widetilde{V}}=\xi-{\widetilde{Z}}$, it follows that
$V-{\widetilde{V}}=w$ and ${\widetilde{V}}$ solves
\begin{equation}\label{e:vtilde}
\dot{\widetilde{V}}+{\widetilde{V}}_{xxxx}+{\widetilde{V}}_{xx}=\alpha{\widetilde{Z}}-{\widetilde{Z}}_{xx}+[({\widetilde{V}}_x+{\widetilde{Z}}_x)^2]_{xx}.
\end{equation}
The corresponding energy functional is given by
\begin{multline*}
\mathcal{E}_t^\alpha(v,z)
= \frac12|v(t)|_{L^2}^2+\\
 +\int_0^t(|v_{xx}|_{L^2}^2-|v_x|_{L^2}^2-\alpha\langle v,z\rangle_{L^2}-\langle v_x,z_x\rangle_{L^2}-\langle 2v_x z_x+(z_x)^2,v_{xx}\rangle_{L^2}),
\end{multline*}
and in particular $\mathcal{E}_t^0=\mathcal{E}_t$.
\begin{lemma}
Let $P$ be an energy martingale solution, then for every $z_0\in H^1_{per}$ and $\alpha\ge0$,
\begin{equation}\label{e:general_ei}
P[\mathcal{E}_t^\alpha({\widetilde{V}},{\widetilde{Z}})\le\mathcal{E}_s^\alpha({\widetilde{V}},{\widetilde{Z}})]=1,
\end{equation}
for almost every $s\ge0$ (including $s=0$) and every $t\ge s$, where ${\widetilde{V}}$,
${\widetilde{Z}}$ have been defined above.
\end{lemma}
\begin{proof}
The proof works as in \cite[Theorem 2.8]{Rom01} and we give just a sketch.
Since ${\widetilde{V}}=V-w$, it follows that
$$
|{\widetilde{V}}(t)|_{L^2}^2=|V(t)|_{L^2}^2+|w(t)|_{L^2}^2-2\langle V(t),w(t)\rangle_{L^2},
$$
and, since by assumptions the energy inequality holds for $V$, it is sufficient to
prove a balance equality for $w$ and $\langle V(t),w(t)\rangle_{L^2}$.
Indeed, it is easy to show by regularisation that
\begin{equation}\label{e:eiforw}
\frac12|w(t)|_{L^2}^2+\int_s^t(|w_{xx}|_{L^2}^2+\alpha\langle{\widetilde{Z}},w\rangle_{L^2})\,ds=\frac12|w(s)|_{L^2}^2
\end{equation}
$P$-a.\ s.\ for all $s\ge0$ and $t\ge s$. We only need to show
that for almost all $s\ge 0$ and $t\ge s$,
\begin{align}\label{e:eiforwv}
&\quad\langle V(t),w(t)\rangle_{L^2}-\langle V(s),w(s)\rangle_{L^2}=\\
&=\!-\!\!\int_s^t\!\!\!\bigl[(2\langle V_{xx},w_{xx}\rangle_{L^2}
          -\langle w_x,V_x+Z_x,\rangle_{L^2}
          +\alpha\langle V,{\widetilde{Z}}\rangle_{L^2}
          -\langle w_{xx},(V_x+Z_x)^2\rangle_{L^2})\bigr].\notag
\end{align}
We sketch the proof of the above formula. Since
$V\in L^\infty(L^2_{per})\cap L^2(H^2_{per})$,
by Lemma \ref{l:tderV} it follows that $\dot V\in L^2(H^{-3}_{per})$.
Moreover, we know that $z_0\in H^1_{per}$ and $\widetilde{Z}_x\in L^4$
and so it is easy to see (by writing the energy balance for $|w_x|_{L^2}^2$)
that $w\in L^2(H^3)$, hence $\dot w\in L^2(H^{-2})$. By slightly adapting
Lemma $1.2$ of Temam \cite[\S III]{Tem}, this implies that
$\langle V,w\rangle_{L^2}$ is differentiable in time with
derivative $\langle\dot V,w\rangle_{H^{-3},H^3}+\langle\dot w, V\rangle_{H^{-2}, H^2}$.
Integration by parts then gives \eqref{e:eiforwv}.

Finally, \textbf{\footnotesize\textsf{[E3]}},
\eqref{e:eiforw} and \eqref{e:eiforwv} together provide \eqref{e:general_ei}.
\end{proof}
\begin{proposition}\label{l:measurable}
Given $z_0\in H^1_{per}$ and $\alpha\ge0$, denote by ${\widetilde{V}}$ and ${\widetilde{Z}}$ the
processes defined above corresponding to $z_0$ and $\alpha$. Then the map
$(t,\omega)\in[0,\infty)\times\Omega\mapsto\mathcal{E}_t^\alpha({\widetilde{V}}(\omega),{\widetilde{Z}}(\omega))$
is progressively measurable and
\begin{itemize}
\item[(i)] for all $0\le s\le t$, the sets
  $E_{s,t}(z_0,\alpha)=\{\mathcal{E}_t^\alpha({\widetilde{V}},{\widetilde{Z}})\le\mathcal{E}_s^\alpha({\widetilde{V}},{\widetilde{Z}})\}$
  are $\mathcal{B}_t$-measurable;
\item[(ii)] for all $t>0$, the sets
  $$E_{t}(z_0,\alpha)=\{\mathcal{E}_t^\alpha({\widetilde{V}},{\widetilde{Z}})\le\mathcal{E}_s^\alpha({\widetilde{V}},{\widetilde{Z}})\textrm{ for a.\ e.\ $s\le t$ (including $0$)}\}$$
  are $\mathcal{B}_t$-measurable;
\item[(iii)] the set
  $$E(z_0,\alpha)=R\cap\{\mathcal{E}_t^\alpha({\widetilde{V}},{\widetilde{Z}})\le\mathcal{E}_s^\alpha({\widetilde{V}},{\widetilde{Z}})\text{ for a.e. }s\ge0\text{ (incl. $0$), all }t\ge s\}$$
  is $\mathcal{B}$-measurable, where
  $$R=\{Z\in L^4_{loc}([0,\infty);W^{1,4}_{per}),\ V\in L^\infty_{loc}([0,\infty);L^2_{per})\cap L^2_{loc}([0,\infty);H^2_{per})\}.$$
\end{itemize}
Moreover, given a \emph{energy martingale solution} $P$, property {\rm\footnotesize\textbf{\textsf{[E3]}}}
is equivalent to each of the following:
\begin{itemize}
\item[\rm\footnotesize\bfseries\textsf{[E3a]}] There are $z_0\in H^1_{per}$ and $\alpha\ge0$ such that
  for each $t>0$ there is a set $T\subset(0,t]$ of null Lebesgue measure and
  $P[E_{s,t}(z_0,\alpha)]=1$ for all $s\not\in T$.
\item[\rm\footnotesize\bfseries\textsf{[E3b]}] There are $z_0\in H^1_{per}$ and $\alpha\ge0$ such that
  for each $t>0$, $P[E_t(z_0,\alpha)]=1$.
\item[\rm\footnotesize\bfseries\textsf{[E3c]}] There are $z_0\in H^1_{per}$ and $\alpha\ge0$ such that $P[E(z_0,\alpha)]=1$.
\end{itemize}
\end{proposition}
\begin{proof}
Measurability of the map $\mathcal{E}^\alpha$ follows from the semi-continuity
properties of the various term of $\mathcal{E}^\alpha$ with respect to the
topology of $\Omega$ (see also Lemma 2.1 of Flandoli \& Romito \cite{FlRo06b}).

The measurability of each $E_{s,t}(z_0,\alpha)$ now follows easily from measurability
of the map $\mathcal{E}^\alpha$. As it regards \emph{(ii)}, fix $t>0$ and notice that
the Borel $\sigma$-algebra of the interval $(0,t)$ is countably generated, so
that if $\mathcal{T}_t$ is a countable basis,
$$
E_t(z_0,\alpha)=E_{0,t}(z_0,\alpha)\cap\bigcap_{T\in\mathcal{T}_t}\{\int_0^t\mathbf{1}_T(s)(\mathcal{E}_t^\alpha({\widetilde{V}},{\widetilde{Z}})-\mathcal{E}_s^\alpha({\widetilde{V}},{\widetilde{Z}}))\,ds\le0\}
$$
and all sets
$\{\int_0^t\mathbf{1}_T(s)(\mathcal{E}_t^\alpha({\widetilde{V}},{\widetilde{Z}})-\mathcal{E}_s^\alpha({\widetilde{V}},{\widetilde{Z}}))\,ds\le0\}$
are $\mathcal{B}_t$-measurable by the measurability of $\mathcal{E}^\alpha$.

We next show \emph{(iii)}. Let $J\subset[0,\infty)$ be a countable dense subset
and define
$$
R_t=\{Z\in L^4_{loc}([0,t);W^{1,4}_{per}),\quad V\in L^\infty(0,t;L^2_{per})\cap L^2(0,t;H^2_{per})\}
$$
(notice that the regularity of $Z$ and $V$ implies that of ${\widetilde{V}}$ and ${\widetilde{Z}}$),
then $R_t\in\mathcal{B}_t$ and, by the lower semi-continuity of the
various terms of $\mathcal{E}_t^\alpha({\widetilde{V}},{\widetilde{Z}})-\mathcal{E}_s^\alpha({\widetilde{V}},{\widetilde{Z}})$
with respect to $t$, it follows that
$$
E(z_0,\alpha)=\bigcap_{t\in J}(R_t\cap E_t(z_0,\alpha))
$$
is $\mathcal{B}$-measurable.
The last statement of the lemma
is now obvious from the above equalities, property
\textbf{\textsf{\footnotesize[E2]}} and Lemma \ref{l:Zreg}.
\end{proof}
\section{Existence of Markov solutions}\label{s:exist}
This section is devoted to the existence of Markov solutions for
equation \eqref{e:dirkeq}. We state the main theorem of this part.
\begin{theorem}\label{t:exmarkov}
There exists a family $(P_x)_{x\in L^2_{per}}$ of probability measures
on $(\Omega,\mathcal{B})$ such that for each $x\in L^2_{per}$, $P_x$
is a energy martingale solution with initial distribution $\delta_x$,
and the a.\ s.\ Markov property holds: there is a set $T_P\subset(0,\infty)$
with null Lebesgue measure such that for all $s\not\in T_P$, all $t\ge s$
and all bounded measurable $\phi:L^2_{per}\to\mathbf{R}$,
$$
\mathbb{E}^P[\phi(\xi_t)|\mathcal{B}_t]
=\mathbb{E}^{P_{\xi_s}}[\phi(\xi_{t-s})].
$$
\end{theorem}
In order to show the theorem, we use the method developed in Flandoli
\& Romito \cite{FlRo06b} (cf. Theorem 2.8). Define for each $x\in L^2_{per}$,
$$
\mathcal{C}(x)=\{\,P\,:\,P\text{ is a energy martingale solution starting at }\delta_x\,\}.
$$
The proof boils down to show that the family $(\mathcal{C}(x))_{x\in L^2_{per}}$
is an \emph{a.\ s.\ pre-Markov family}. We recall here the
various properties that we need to show to prove the statement
(see also Definition 2.5 of Flandoli \& Romito \cite{FlRo06b}).
\begin{enumerate}
\item Each $\mathcal{C}(x)$ is non-empty, compact and convex, and
  the map $x\to\mathcal{C}(x)$ is measurable with respect to the
  Borel $\sigma$-fields of the space of compact subsets of $\Pr(\Omega)$
  (endowed with the Hausdorff measure).
\item For each $x\in L^2_{per}$ and all $P\in\mathcal{C}(x)$,
  $P[C([0,\infty);L^2_{per,\sigma})]=1$, where $L^2_{per,\sigma}$
  is the space $L^2_{per}$ with the weak topology.
\item For each $x\in L^2_{per}$ and $P\in\mathcal{C}(x)$ there is a set
  $T\subset(0,\infty)$ with null Lebesgue measure, such that for all
  $t\not\in T$ the following properties hold:
  \begin{enumerate}
  \item\emph{(disintegration)} there exists $N\in\mathcal{B}_t$ with $P(N)=0$
    such that for all $\omega\not\in N$
    \begin{center}
    $\omega(t)\in L^2_{per}$\qquad and \qquad $P|_{\mathcal{B}_t}^\omega\in\Phi_t\mathcal{C}(\omega(t))$;
    \end{center}
  \item\emph{(reconstruction)} for each $\mathcal{B}_t$-measurable map
    $\omega\mapsto Q_\omega:\Omega\to\Pr(\Omega^t)$ such that there
    is $N\in\mathcal{B}_t$ with $P(N)=0$ and for all $\omega\not\in N$,
    \begin{center}
    $\omega(t)\in L^2_{per}$\qquad and \qquad $Q_\omega\in\Phi_t\mathcal{C}(\omega(t))$;
    \end{center}
    then $P\otimes_t Q\in\mathcal{C}(x)$.
  \end{enumerate}
\end{enumerate}
The validity of this statement is verified in the following lemmas.
\begin{lemma}[Continuity lemma]
\label{lem:cont}
For each $x\in L^2_{per}$, the set $\mathcal{C}(x)$ is non-empty,
convex and for all $P\in\mathcal{C}(x)$,
$$
P[C([0,\infty);L^2_{per,\sigma}]=1.
$$
\end{lemma}
\begin{proof}
Existence of weak martingale mild solutions is proved in Bl\"omker \& Gugg
\cite{BlGu}, using standard spectral Galerkin methods. This is similar to Lemma \ref{lem:comp}.

By Remark \ref{r:equiv}, this implies existence of weak
martingale solutions according to Definition \ref{d:weakms}. In order
to prove the energy inequality of Definition \ref{d:energyms}, one can
proceed as in the next lemma (where it is proved in a slightly more
general situation).

Next, it is easy to show that $\mathcal{C}(x)$ is convex, since
all requirements of both Definitions \ref{d:weakms} and \ref{d:energyms}
are linear with respect to measures $P\in\mathcal{C}(x)$. Finally,
if $P\in\mathcal{C}(x)$, we know by Lemma \ref{l:Zreg} that,
under $P$, the process $Z$ is weakly continuous.
Moreover, by property \textbf{\textsf{\footnotesize[E2]}} of Definition
\ref{d:energyms}, $V$ is also weakly continuous
and in conclusion $C([0,\infty);L^2_{per,\sigma})$ is a full set.
\end{proof}
\begin{lemma}[Compactness lemma]\label{lem:comp}
For each $x\in L^2_{per}$, the set $\mathcal{C}(x)$ is compact and the
map $x\mapsto\mathcal{C}(x)$ is Borel measurable. 
\end{lemma}
\begin{proof}
Following Lemma 12.1.8 of Stroock \& Varadhan \cite{StVa}, it is sufficient
to prove that for each sequence $(x_n)_{n\in\mathbf{N}}$ converging to $x$
in $L^2$ and for each $P_n\in\mathcal{C}(x_n)$, the sequence
$(P_n)_{n\in\mathbf{N}}$ has a limit point $P$, with respect to weak
convergence of measures, in $\mathcal{C}(x)$.

Let $x_n\to x$ in $L^2_{per}$ and let $P_n\in\mathcal{C}(x_n)$. By
Theorem \ref{t:tightness}, $(P_n)_{n\in\mathbf{N}}$ is tight on
$\Omega\cap L^2_{loc}([0,\infty);H^1_{per})$. Hence, up to a
sub-sequence that we keep denoting by $(P_n)_{n\in\mathbf{N}}$,
it follows that $P_n\rightharpoonup P$, for some $P$. It remains
to show that $P\in\mathcal{C}(x)$. Therefore, we verify all properties of 
Definitions \ref{d:weakms} and \ref{d:energyms}.

We start by proving \textbf{\textsf{\footnotesize[W2]}} for $P$. Given
$\varphi\in\mathcal{D}^\infty$, we know that for each $n\in\mathbf{N}$
the process $(|\mathcal{Q}^{1/2}\varphi|_{L^2}^{-1}M_t^\varphi,\mathcal{B}_t,P_n)_{t\ge0}$
is a one-dimensional standard Brownian motion. Now, since
$P_n\rightharpoonup P$ and $M_t^\varphi$ is continuous
(with respect to both $\omega$ and $t$), it follows that
the law of $M^\varphi$ under $P$ is that of a
one-dimensional standard Brownian motion.

Property \textbf{\textsf{\footnotesize[W3]}} is obvious, since
the marginals of $P_n$ at time $0$ converge, by assumption, to
both $\delta_x$ and the marginal of $P$ at time $0$, hence they
coincide and $P$ is started at $\delta_x$.

In order to prove the other properties, 
we rely on tightness from Theorem \ref{t:tightness} with
$K=\log(1+|x|_{L^2}^2)^\kappa$, and use the classical
Skorokhod theorem: there exist a probability space
$(\Sigma,\mathcal{F},\mathbb{P})$ and random variables
$(h^{(n)},z^{(n)})_{n\in\mathbf{N}}$ and
$(h^{(\infty)},z^{(\infty)})_{n\in\mathbf{N}}$
such that each $(h^{(n)},z^{(n)})$ has the same law
of $(\xi,Z)$ under $P_n$ (and similarly
for $(h^{(\infty)},z^{(\infty)})$) and $h^{(n)}\to h^{(\infty)}$
in $\Omega\cap L^2_{loc}([0,\infty);H^1)$ and
$z^{(n)}\to z^{(\infty)}$ in $L^{\frac{16}3}(0,T;W^{1,4})$,
$\mathbb{P}$-a.\ s.. In particular, $v^{(n)}=h^{(n)}-z^{(n)}$ has
the same law of $V$ under $P_n$ (and so is for
$v^{(\infty)}=h^{(\infty)}-z^{(\infty)}$ and $V$
under $P$).

In order to prove \textbf{\textsf{\footnotesize[W1]}}, it
is sufficient to show that 
$$
\mathbb{P}[\|h\|_{L^2(0,T;H^1)}>K]\to0
\qquad\text{as $K\uparrow\infty$ for all $T>0$.}
$$
By \eqref{e:logtightforh},
we know that $\mathbb{E}^\mathbb{P}[\log\bigl(1+\int_0^T|h^{(n)}|^2_{H^1}\bigr)]\le C_T$,
so that Fatou's lemma implies a similar estimate for $h^{(\infty)}$
and Chebychev inequality gives the result.

One can proceed similarly to prove \textbf{\textsf{\footnotesize[E2]}},
using \eqref{e:logtightforv} and the fact that norms
in $L^\infty(0,T;L^2)$ and in $L^2(0,T;H^2)$ are lower semi-continuous
with respect to the topology where $v^{(n)}\to v^{(\infty)}$.

In order to prove \textbf{\textsf{\footnotesize[E3]}}, we show
that property \textbf{\textsf{\footnotesize[E3a]}} (with $z_0=0$ and $\alpha=0$)
of Proposition \ref{l:measurable} holds true. Fix $t>0$. A first useful fact is
that $v^{(n)}$ converges weakly in $L^2(0,t;H^2)$ to $v^{(\infty)}$. Indeed, we
can use \textbf{\textsf{\footnotesize[E3]}}, applied to each $v^{(n)}$, and the
convergence of $z^{(n)}$ to show that $(v^{(n)})_{n\in\mathbf{N}}$ is bounded in
$L^2(0,t;H^2)$, $\mathbb{P}$-a.\ s. (the bound follows from an inequality for each
$v^{(n)}$ which can be obtained from the energy inequality as \eqref{e:vbound}
in Lemma \ref{l:vbound}). It follows then that $v^{(n)}\rightharpoonup v^{(\infty)}$,
in $L^2(0,t;H^2)$, since we already know that $v^{(n)}$ converges to $v^{(\infty)}$
in $L^2(0,t;H^1_{per})$.

A second useful fact is that there is a null Lebesgue set $S\subset (0,t]$ such
that for all $s\not\in S$,
\begin{equation}\label{e:claimE3}
\mathbb{P}\left[|v^{(n')}(s)|_{L^2}\to|v^{(\infty)}(s)|_{L^2}
\text{ for a subsequence }v^{(n')}\right]=1.
\end{equation}
Note that this does not imply a.s.\ convergence for a subsequence,
as the subsequence may depend on  $\sigma\in\Sigma$.

To prove \eqref{e:claimE3} note that $v^{(n)}\to v^{(\infty)}$, $\mathbb{P}$-a.\ s.\ in
$L^2(0,t;L^2_{per})$, and so
$$
\mathbb{E}^\mathbb{P}[\log(1+\frac1t\int_0^t|v^{(n)}-v^{(\infty)}|^2_{L^2}\,ds)]\to0.
$$
This follows from uniform
bounds on higher moments from \eqref{e:logtightforv} with $\kappa>1$. 
By Jensen inequality,
$$
\mathbb{E}^\mathbb{P}[\frac1t\int_0^t\log(1+|v^{(n)}-v^{(\infty)}|^2_{L^2})\,ds]
\le \mathbb{E}^\mathbb{P}[\log(1+\frac1t\int_0^t|v^{(n)}-v^{(\infty)}|^2_{L^2}\,ds)],
$$
and so there are a set $S\subset(0,t]$ (notice that $0\not\in S$ since we already
know that $v^{(n)}(0)\to v^{(\infty)}(0)$) and a subsequence $v^{(n')}$
such that 
$$
\mathbb{E}^\mathbb{P}[\log(1+|v^{(n')}(s)-v^{(\infty)}(s)|^2_{L^2})]\to0
\qquad\text{for all }s\not\in S.
$$
From this claim \eqref{e:claimE3} now easily follows, possibly by taking a
further sub-sequence depending on $\sigma\in\Sigma$.

We are now able to prove \textbf{\textsf{\footnotesize[E3a]}} for $P$
(with $z_0=0$ and $\alpha=0$). We know that for each $n\in\mathbf{N}$
there is a null Lebesgue set $T_n\subset(0,t]$ such that
$\mathbb{P}[\mathcal{E}_t(v^{(n)},z^{(n)})\le\mathcal{E}_s(v^{(n)},z^{(n)})]=1$,
for all $s\not\in T_n$. Let $T=S\cup\bigcup T_n$ and consider $s\not\in T$. Since
we know that $\mathcal{E}_t(v^{(n)},z^{(n)})\le\mathcal{E}_s(v^{(n)},z^{(n)})$ holds
$\mathbb{P}$-a.\ s.\ for all $n\in\mathbf{N}$, by passing to the limit $n\to\infty$
and using all the convergence information we have collected, it follows that
$\mathbb{P}[\mathcal{E}_t(v^{(\infty)},z^{(\infty)})\le\mathcal{E}_s(v^{(\infty)},z^{(\infty)})]=1$.
\end{proof}
Before stating the next two lemmas (which contain the multi-valued
form of the Markov property), we need to analyse what happens to
processes $W$, $Z$ and $V$ under the action of the forward shift
$\Phi_u$, for a given $u$. First, given $s\ge0$ and $z_0\in H^1_{per}$,
denote by $Z(t,\cdot|s,z_0)$ the Ornstein-Uhlenbeck process starting
in $z_0$ at  time $s$, namely
$$
Z(t,\cdot|s,z_0)=\mathrm{e}^{A(t-s)}z_0+\sum\alpha_k\int_s^t\mathrm{e}^{(t-r)\lambda_k}d\beta_k(r)e_k.
$$
In particular, we have that $Z(t,\cdot|0,0)=Z(t,\cdot)$. Set
moreover $V(t,\cdot|s,z_0)=\xi-Z(t,\cdot|s,z_0)$. Now, from
\textbf{\textsf{\footnotesize[W2]}} and \eqref{e:defW}
it is easy to verify that, for all $\omega\in\Omega^u$,
$$
W(t,\Phi_u^{-1}(\omega))=W(t+u,\omega)-W(u,\omega),
$$
and it depends only on the values of $\omega$ in $[u,u+t]$.
Similarly,
\begin{align}\label{e:shifted}
Z(\Phi_u^{-1}(\omega),t|s,z_0)=Z(\omega,t+u,|s+u,z_0),\notag\\
V(\Phi_u^{-1}(\omega),t|s,z_0)=V(\omega,t+u|s+u,z_0).
\end{align}
\begin{lemma}[Disintegration lemma]
For every $x\in L^2_{per}$ and $P\in\mathcal{C}(x)$, there is a set $T\subset(0,\infty)$,
with null Lebesgue measure, such that for all $t\not\in T$ there is $N\in\mathcal{B}_t$,
with $P[N]=0$, such that for all $\omega\not\in N$,
$$
\omega(t)\in L^2_{per}
\qquad\text{and}\qquad
P|_{\mathcal{B}_t}^\omega\in\Phi_t\mathcal{C}(\omega(t)).
$$
\end{lemma}
\begin{proof}
Fix $x\in L^2_{per}$ and $P\in\mathcal{C}(x)$, let $T_P$ be the set of
exceptional times of $P$, as given by \textbf{\footnotesize\textsf{[E3]}} of Definition
\ref{d:energyms}, and fix $u\not\in T_P$. Let
$(P|^\omega_{\mathcal{B}_u})_{\omega\in\Omega}$ be a regular conditional
probability distribution of $P$ given $\mathcal{B}_u$. We aim to show
that there is a $P$-null set $N\in\mathcal{B}_u$ such that
$\omega(u)\in L^2_{per}$ and
$P|^\omega_{\mathcal{B}_u}\in\Phi_u\mathcal{C}(\omega(u))$
for all $\omega\not\in N$. We shall find
$N=N_{\textsf{\tiny[E1]}}\cup N_{\textsf{\tiny[E2]}}\cup N_{\textsf{\tiny[E3]}}$,
corresponding to \emph{bad} sets for each property.

We only prove \textbf{\footnotesize\textsf{[E2]}} and
\textbf{\footnotesize\textsf{[E3]}} for the conditional distributions $P|^\omega_{\mathcal{B}_u}$,
the proof of the other properties being entirely similar to Lemma
4.4 of Flandoli \& Romito \cite{FlRo06b}. We start by
\textbf{\footnotesize\textsf{[E2]}}. We need to show that
$P|^\omega_{\mathcal{B}_u}[V(\cdot,\Phi_u^{-1}(\cdot))\in S_{[0,\infty)}]=1$
or, equivalently, by \eqref{e:shifted}, that 
$P|^\omega_{\mathcal{B}_u}[V(\cdot,\cdot|u,0)\in S_{[u,\infty)}]=1$,
where we have set, for brevity, $S_J=L^\infty_{loc}(J;L^2_{per})\cap L^2_{loc}(J;H^2_{per})$,
for any interval $J\subset[0,\infty)$.
Set
\begin{align}\label{e:Vdisreg}
\mathcal{S}_u&=\{V\in S_{[0,u]}\text{ and }\mathrm{e}^{At}Z(u,\cdot)\in S_{[0,\infty)}\},\notag\\
\mathcal{S}^u&=\{V(\cdot,\cdot|u,0)\in S_{[u,\infty)}\},
\end{align}
then $\mathcal{S}_u\in\mathcal{B}_u$ and $\mathcal{S}^u\in\mathcal{B}^u$, since by definition
$V$ and $Z$ are adapted. Moreover, since
$V(t+u,\omega)=V(t+u,\omega|u,0)-\mathrm{e}^{At}Z(u,\omega)$,
it follows from \textbf{\textsf{\footnotesize[E2]}} for $P$,
Lemma \ref{lem:ZregLinfty} and the regularity properties of
the semigroup $\mathrm{e}^{At}$, that $\mathcal{S}_u\cap\mathcal{S}^u$
is a $P$-full set and so, by disintegration,
$$
1=P[\mathcal{S}_u\cap\mathcal{S}^u]=\int_{\mathcal{S}_u}P|_{\mathcal{B}_u}^\omega[\mathcal{S}^u]\,P(d\omega).
$$
Thus, there is a $P$-null set $N_{\textsf{\tiny[E2]}}\in\mathcal{B}_u$ such that
$P|_{\mathcal{B}_u}^\omega[\mathcal{S}^u]=1$ for all $\omega\not\in N_{\textsf{\tiny[E2]}}$.

We finally prove \textbf{\textsf{\footnotesize[E3c]}} (cf. Proposition \ref{l:measurable})
for the conditional probabilities. Set
\begin{align*}
A  &=\{\mathcal{E}_t(V,Z)\le\mathcal{E}_s(V,Z)\text{\small\ for a.e.\ }s\ge0\text{\small\ (including $0$, $u$), all }t\ge s\}\\
A_u&=\{\mathcal{E}_t(V,Z)\le\mathcal{E}_s(V,Z)\text{\small\ for a.e.\ }s\in[0,u]\text{\small\ (incl. $0$, $u$), all }t\in[s,u]\},
\end{align*}
where, for the sake of simplicity, in the definitions of the above sets we have omitted
the information on regularity for $V$ and $Z$, which are essential to ensure measurability
(compare with Proposition \ref{l:measurable}).
They can be treated as in the proof of \textbf{\textsf{\footnotesize[E2]}} above. We have
$A_u\in\mathcal{B}_u$ and $P[A]=P[A_u]=1$, since $u\not\in T_P$.
Now, if $\overline{\omega}\in A_u\cap\{Z\in H^1_{per}\}$
(which is again a $P$ full set by Lemma \ref{lem:ZregLinfty}), set
$$
B(\overline{\omega})=A\cap\{\omega:\omega=\overline{\omega}\text{ on }[0,u]\}
$$
and notice that, for such $\overline{\omega}$, $B(\overline{\omega})$ is equal to
$$
\{\text{\small$\mathcal{E}_t(V_{\overline{\omega}},Z_{\overline{\omega}})
            \le\mathcal{E}_s(V_{\overline{\omega}},Z_{\overline{\omega}})$
              for a.\ e.\ $s\ge u$ (incl. $u$), all $t\ge s$}\}
$$
since $V(t+u,\omega)=V(t+u,\omega|u,Z(u,\omega))$ (a similar relation holds for $Z$ as well),
and we have set $V_{\overline{\omega}}(\cdot)=V(\cdot|u,Z(u,\overline{\omega}))$
and $Z_{\overline{\omega}}(\cdot)=Z(\cdot|u,Z(u,\overline{\omega}))$. Moreover, the map
$$
\omega\to\mathbf{1}_{A_u\cap\{Z\in H^1_{per}\}}(\omega)P|_{\mathcal{B}_u}^\omega[B(\omega)]
$$
is $\mathcal{B}_u$-measurable, since $P|_{\mathcal{B}_u}^\omega[B(\omega)]=P|_{\mathcal{B}_u}^\omega[A]$
for all $\omega\in A_u\cap\{Z\in H^1_{per}\}$. Now, by \textbf{\textsf{\footnotesize[E3c]}}
for $P$ (with $z_0=0$ and $\alpha=0$) and disintegration,
$$
1=P[A]
=\mathbb{E}^P[\mathbf{1}_{A_u\cap\{Z\in H^1_{per}\}}(\cdot)P|_{\mathcal{B}_u}^\cdot[B(\cdot)]],
$$
and so there is $N_{\textsf{\tiny[E3]}}\in\mathcal{B}_u$ such that
$P|_{\mathcal{B}_u}^\omega[B(\omega)]=1$ for all $\omega\not\in N_{\textsf{\tiny[E3]}}$
or, in different words, such that \textbf{\textsf{\footnotesize[E3c]}} holds
(with $z_0=Z(u,\omega)$ and $\alpha=0$) for $P|_{\mathcal{B}_u}^\omega$
for all $\omega\not\in N_{\textsf{\tiny[E3]}}$.
\end{proof}
\begin{lemma}[Reconstruction lemma]
For every $x\in L^2_{per}$ and $P\in\mathcal{C}(x)$, there is a set $T\subset(0,\infty)$,
with null Lebesgue measure, such that for each $t\not\in T$, for each
$\mathcal{B}_t$-measurable map $\omega\mapsto Q_\omega:\Omega\to\Pr(\Omega^t)$
such that there is $N\in\mathcal{B}_t$, with $P[N]=0$, and for all $\omega\not\in N$,
$$
\omega(t)\in L^2_{per}
\qquad\text{and}\qquad
Q_\omega\in\Phi_t\mathcal{C}(\omega(t)),
$$
then $P\otimes_t Q\in\mathcal{C}(x)$.
\end{lemma}
\begin{proof}
Let $x\in L^2_{per}$, $P\in\mathcal{C}(x)$, $T_P$ be the set of exceptional
times of $P$ and fix $u\not\in T_P$. Let $(Q_\omega)_{\omega\in\Omega}$
be a $\mathcal{B}_u$-measurable map and $N_Q$ a $P$-null set such that
$\omega(u)\in L^2_{per}$ and $Q_\omega\in\Phi_u\mathcal{C}(\omega(u))$
for all $\omega\not\in N_Q$. In order to verify that $P\otimes_u Q\in\mathcal{C}(x)$,
we only check properties \textbf{\textsf{\footnotesize[E2]}} and
\textbf{\textsf{\footnotesize[E3]}}, since the proof of
\textbf{\textsf{\footnotesize[E1]}} can be carried on as in Flandoli
\& Romito \cite[Lemma 4.5]{FlRo06b}.

We start by \textbf{\textsf{\footnotesize{[E2]}}}.
Consider again sets $\mathcal{S}_u\in\mathcal{B}_u$ and $\mathcal{S}^u\in\mathcal{B}^u$
defined in \eqref{e:Vdisreg} and notice that, by \textbf{\textsf{\footnotesize{[E2]}}}
for $Q_\omega$, for each $\omega\not\in N_Q$ we have that $Q_\omega[\mathcal{S}^u]=1$.
Moreover, by \textbf{\textsf{\footnotesize{[E2]}}} for $P$, Lemma \ref{lem:ZregLinfty}
and the regularity properties of the semigroup $\mathrm{e}^{At}$, it follows that
$P[\mathcal{S}_u]=1$. Finally, since we know that
$V(t+u,\omega)=V(t+u,\omega|u,0)-\mathrm{e}^{At}Z(u,\omega)$,
it follows easily that $\mathcal{S}_u\cap\mathcal{S}^u=\{V\in S_{[0,\infty)}\}$
and so
$$
(P\otimes_u Q)[V\in S_{[0,\infty)}]
=(P\otimes_u Q)[\mathcal{S}_u\cap\mathcal{S}^u]
=\int_{\mathcal{S}_u}Q_\omega[\mathcal{S}^u]\,P(d\omega)
=1.
$$

We next prove \textbf{\textsf{\footnotesize{[E3]}}}. Again, we prove
it by means of \textbf{\textsf{\footnotesize{[E3c]}}}, thanks to Proposition
\ref{l:measurable}. Define $A$ and $A_u$ as in the proof of the previous lemma
(the regularity conditions on $Z$ and $V$ are again omitted). Since $u\not\in T_P$
and $A_u\in\mathcal{B}_u$, we know that $(P\otimes_u Q)[A_u]=P[A_u]=1$. Moreover, by Lemma
\ref{lem:ZregLinfty}, there is a $P$-null set $N\in\mathcal{B}_u$ such that
$Z(u,\omega)\in H^1_{per}$ for all $\omega\not\in N$. For each $\overline{\omega}\not\in N$,
define $B(\overline{\omega})=A\cap\{\omega:\omega=\overline{\omega}\text{ on }[0,u]\}$
and notice that, if $\omega\in A_u\cap(N\cap N_Q)^c$ (which is again a
$\mathcal{B}_u$-measurable $(P\otimes_u Q)$-full set), then by \textbf{\textsf{\footnotesize{[E3c]}}}
(with $z_0=Z(u,\omega)$ and $\alpha=0$) for $Q_\omega$ it follows that $Q_\omega[B(\omega)]=1$.
The map $\omega\mapsto\mathbf{1}_{A_u\cap(N\cap N_Q)^c}(\omega)Q_\omega[B(\omega)]$
is then trivially $\mathcal{B}_u$-measurable and equal to $1$, $P$-a.\ s..
Moreover, we have that $Q_\omega[A]=Q_\omega[B(\omega)]=1$ for all
$\omega\in A_u\cap(N\cap N_Q)^c$ and so
$$
(P\otimes_u Q)[A]
=\mathbb{E}^P\bigl[\mathbf{1}_{A_u\cap(N\cap N_Q)^c}Q_\cdot[B(\cdot)]\bigr]
=P[A_u\cap(N\cap N_Q)^c]
=1.
$$
In conclusion, \textbf{\textsf{\footnotesize{[E3c]}}} (with $z_0=0$ and $\alpha=0$)
holds true for $P\otimes_u Q$.
\end{proof}
\section{The strong Feller property}\label{s:sfeller}
Throughout this section we shall assume that the noise is non-degenerate.
This is summarised by the following assumption.
\begin{assumption}\label{a:noiseass}
The operator $\mathcal{Q}^{-\frac12}$ is bounded, where $\mathcal{Q}$
is the covariance of the noise. In different words,
$$
\alpha_k\ge\delta>0,
$$
for some constant $\delta$, where $\alpha_k^2$ are the eigenvalues of $\mathcal{Q}$.
\end{assumption}
\begin{theorem}\label{t:strongFeller}
Under the above assumption, any a.\ s.\ Markov family $(P_x)_{x\in L^2_{per}}$
of energy martingale solutions defines a Markov semigroup that has the
$H^1$-strong Feller property.
\end{theorem}
\begin{proof}[Proof of Theorem \ref{t:strongFeller}]
We mainly rely on \cite{FlRo06b} and \cite{FlRo06c}. 
Let $(P_x)_{x\in L^2_{per}}$ be an a.\ s.\ Markov family of energy martingale
solution and denote by $(\mathcal{P}_t)_{t\ge0}$ the corresponding (a.\ s.) semigroup
generated by  $P_x$. Then the claim follows from the following lemma.
\begin{lemma}\label{lem:SF-cruc}
There is an $\epsilon = \epsilon(|h|_{H^1},R)\to 0$ for $h\to 0$
such that
\begin{equation}\label{e:SF-cruc}
|\mathcal{P}_\epsilon\varphi(x+h)-\mathcal{P}_\epsilon\varphi(x)| 
\le C |h|_{H^1} \log(1/|h|_{H^1})
\end{equation}
for all $ |h|_{H^1}\le 1$, all $\varphi\in L^\infty(H^1)$ with $|\varphi|_{L^\infty}\le 1$, 
and all $|x|_{H^1}\le R/4$ for some sufficiently large $R$.
\end{lemma}
With this lemma at hand, we define for $\varphi\in L^\infty(H^1)$ with $|\varphi|_{L^\infty}=1$
and $h$ (i.e., $\epsilon$) sufficiently small $\varphi_*=\mathcal{P}_{t-\epsilon}\varphi\in L^\infty(H^1)$
with $|\varphi_*|_{L^\infty}\le 1$. Thus
\begin{equation}
|\mathcal{P}_t\varphi(x+h)-\mathcal{P}_t \varphi(x)| 
\le  |\mathcal{P}_\epsilon\varphi_*(x+h)-\mathcal{P}_\epsilon \varphi_*(x)| \\
\le  C |h|_{H^1} \log(1/|h|_{H^1}).
\end{equation}
This implies strong Feller for $\mathcal{P}_t$.
\end{proof}
Following the arguments of \cite{FlRo06b} and \cite{FlRo06c}
it is enough to prove strong Feller for 
the following \emph{regularised} problem
\begin{equation}\label{e:reg}
\partial_t u=-u_{xxxx}+(-u+(u_x)^2)_{xx}\chi_\rho(|u|^2_{H^1})+\partial_t W
\end{equation}
where $\chi_\rho\in C^\infty$ is a cut-off function such that 
$\chi_\rho\equiv1$ on $[0,\rho^2]$ and $\chi_\rho\equiv0$ on $[2\rho^2,\infty)$.
For all $\zeta\ge0$ we have
$$ 
|\chi_\rho(\zeta)| \le 1, \quad 
|\chi_\rho'(\zeta)| \le C \rho^{-2},\quad 
|\chi_\rho(\zeta^2)\zeta^p| \le C \rho^{p},\quad 
|\chi_\rho'(\zeta^2)\zeta^p| \le C \rho^{p-2}.
$$
Let $P^{(\rho)}_x$ be the (unique) Markov energy martingale solution solution
of the regularised problem \eqref{e:reg}. This is well defined, as we can solve
\eqref{e:reg} path-wise. The mild solution of \eqref{e:reg} is given by
\begin{equation}\label{e:regmild}
u(t)=\mathrm{e}^{tA}u(0) -\int_0^t\partial_x^2\mathrm{e}^{(t-\tau)A}F(u(\tau)) d\tau + Z(t)
\end{equation}
where $Z$ has been defined in \eqref{e:defZ} and
$$
F(u)= (-u+(u_x)^2)\chi_\rho(|u|^2_{H^1}).
$$
Using the embedding of $L^1$ into $H^{-1+4\gamma}$ for $\gamma\in(0,\frac18)$,
we can easily check that 
\begin{equation}\label{e:bouF}
|F(u)-F(v)|_{H^{-1+4\gamma}}  \le C_\rho |u-v|_{H^1} 
\quad\mathrm{and}\quad
|F(u)|_{H^{-1+4\gamma}}
\le C(\rho+\rho^2).
\end{equation}
Now uniqueness for \eqref{e:reg} in $C^0([0,\infty),H^1_{per})$ 
follows from standard path-wise fixed point arguments. 
The proof is straightforward as we can rely on one hand on
$F$ being Lipschitz and bounded, and on the other hand $\mathrm{e}^{tA}$ 
generates an analytic semigroup such that
$$ 
|\mathrm{e}^{tA}w|_{H^1} \le |w|_{H^1}
\quad\mathrm{and}\quad
|\partial_x^2\mathrm{e}^{tA}w|_{H^1} \le M (1+t^{\gamma-1}) |w|_{H^{-1+4\gamma}}
$$
(see for example Henry \cite{henry} or  Pazy \cite{pazy}, Lunardi \cite{lunardi}).

Next, define 
\begin{equation}\label{e:blowuptime}
\tau_\rho
=\inf\{t>0:\text{\small the solution of \eqref{e:reg} is bounded in $H^1$ on $[0,t]$ by $\rho$}\}
\end{equation}
Thus the solution of the regularised problem coincides with the energy solution
up to $\tau_\rho$ and in view of \eqref{e:SF-cruc} we have
\begin{equation}\label{e:SF-transreg}
|\mathcal{P}_\epsilon\varphi(x+h)-\mathcal{P}_\epsilon\varphi(x)| 
\le 2\big(P_x[\tau_\rho<\epsilon]+ P_{x+h}[\tau_\rho<\epsilon]\big)
+|\mathcal{P}^{(\rho)}_\epsilon\varphi(x+h)-\mathcal{P}^{(\rho)}_\epsilon\varphi(x)|,
\end{equation}
where $\mathcal{P}^{(\rho)}$ is the semigroup generated by 
\eqref{e:reg} or \eqref{e:regmild}, respectively.

In order to prove Lemma \ref{lem:SF-cruc} we need the 
following two lemmas.
\begin{lemma}\label{lem:SF-SFreg}
There is a $p>1$ sufficiently large,
such that for $\rho\ge 1$ and $t\le1$ 
$$
|\mathcal{P}^{(\rho)}_t\varphi(x+h)-\mathcal{P}^{(\rho)}_t\varphi(x)|
\le \frac{C}t |h|_{H^{-1}}  e^{ct\rho^p} 
$$ 
for all $x,h\in H^1$.
\end{lemma}
\begin{lemma}\label{lem:SF-stop}
There is a small constant $c_\tau$ depending on $\gamma$, and $M$ such that
for all $\rho\ge 1$, $\epsilon\in(0,1]$, $u_0$ such that $|u_0|_{H^1}\le \rho/4+1$,
we have
$$
P_{u_0}[\tau_\rho\ge \epsilon] 
\ge P_{u_0}\big[\sup_{t\in[0,\epsilon]} |Z(t)|_{H^1} \le \rho/4 \big]
$$
for all $\epsilon \le C_\tau \rho^{-2/\gamma}$.
\end{lemma}
Using arguments analogous to \cite[Prop.15]{FlRo06c}
we immediately obtain
\begin{cor}\label{lem:SF-stopcor}
There are two constant $c$, $C>0$ depending on $\gamma$ and $M$ such that
for all $\rho\ge1$, $\epsilon\in(0,1]$, $u_0$ such that $|u_0|_{H^1}\le \rho/4+1$,
we have
$$
P_{u_0}[\tau_\rho\ge \epsilon] \le C\mathrm{e}^{-c\rho^2/\epsilon}
$$
for all $\epsilon \le c_\tau \rho^{-2/\gamma}$.
\end{cor}
\begin{proof}[Proof of Lemma \ref{lem:SF-cruc}]
For $h,x\in H^1$ such that $|x|_{H^1}\le \rho/4$ and   $|h|_{H^1}\le 1$,
we can apply Corollary \ref{lem:SF-stopcor} for $u_0=x$ and $u_0=x+h$.
From \eqref{e:SF-transreg} together with Lemma \ref{lem:SF-SFreg} and the
embedding of $H^1$ into $H^{-1}$
for $\epsilon \le \min\{1,c_\tau \rho^{-2/\gamma}\}$, $\rho\ge \max\{4|x|_{H^1},1\}$,
$t\le 1$,
\begin{equation}
|\mathcal{P}_\epsilon\varphi(x+h)-\mathcal{P}_\epsilon\varphi(x)| 
\le C\mathrm{e}^{-c\rho^2/\epsilon}+C |h|_{H^1} \frac1t\mathrm{e}^{ct\rho^p}.
\end{equation}
Thus, if we fix for a suitable constant $C>0$
$$
\epsilon= \min\Big\{1;\ \frac{C}{\rho^q \ln(1/|h|_{H^1})}\Big\}
\qquad \mathrm{for\ some}\ 
q > \max\{p,2/\gamma\},
$$
then we obtain 
$$|\mathcal{P}_\epsilon\varphi(x+h)-\mathcal{P}_\epsilon\varphi(x)| 
\le C |h|_{H^1} \ln(1/|h|_{H^1}).
$$
\end{proof}
The remainder of the section is devoted to the proof of the two remaining lemmas.
\begin{proof}[Proof of Lemma \ref{lem:SF-stop}]
First from \eqref{e:regmild} for $t\le1$
$$ 
|u(t)|_{H^1} \le |u(0)|_{H^1}+C \int_0^t (t-s)^{\gamma-1} |F(u)|_{H^{-1+4\gamma}}\,ds + |Z(t)|_{H^1}
$$
Thus from \eqref{e:bouF} for $t\le \min\{1,\tau_\rho\}$ and $\rho\ge1$
$$ 
|u(t)|_{H^1} \le \rho/4 + C \tau_\rho^\gamma \rho^2 +|Z(t)|_{H^1}
$$
which easily implies the claim.
\end{proof}


\begin{proof}[Proof of Lemma \ref{lem:SF-SFreg}]
We proceed analogous to the proof of \cite[Proposition 5.12]{FlRo06b}.
For every $v\in H^1_{per}$, let $u(t,v)$ be the
solution to equations \eqref{e:reg} with $u(0,v)=v$. By
the Bismut, Elworthy \& Li formula,
$$
D_y(\mathcal{P}_t^{(\rho)}\varphi)(v)
=\frac1{t}\mathbb{E}[\varphi(u(t,v))\int_0^t\langle\mathcal{Q}^{-1}D_y u(s,v),dW(s)\rangle_{L^2}]
$$
Now Burkholder, Davis \& Gundy inequality states 
$$ \mathbb{E} \sup_{t\in[0,T]}\Big| \int_0^t \langle f(s),dW(s) \rangle_{L^2}\Big|^p
\le C \mathbb{E} \Big(\int_0^T |\mathcal{Q}^{1/2}f(t)|^2_{L^2} dt\Big)^{p/2}
$$
and thus, for $|\varphi|_\infty\leq1$,
\begin{equation}\label{e:SFso}
|(\mathcal{P}_t^{(\rho)}\varphi)(v+h)-(\mathcal{P}_t^{(\rho)}\varphi)(v)|
\le\frac{C}{t}\sup_{\eta\in[0,1]}
\mathbb{E}\Big[\Big(\int_0^t|\mathcal{Q}^{-\frac12} D_h u(s,v+\eta h)|_{L^2}^2\,ds\Big)^{\frac12}\Big].
\end{equation}
Now $\psi(t)=D_h u(t,v+\eta h)$ with $\psi(0)=\eta h$ solves
\begin{equation}\label{e:refdiff}
\partial_t \psi=-\psi_{xxxx}+\partial_x^2 DF(u)[\psi]
\end{equation}
with
$$
DF(u)[\psi]=
-(\psi+2u_x\psi_x)\chi_\rho(|u|^2_{H^1})
-2(u+(u_x)^2)\chi_\rho'(|u|^2_{H^1})\langle u,\psi\rangle_{H^1}.
$$
The following arguments are only formal, but as we are working with unique solutions
they can all be made rigorous by Galerkin approximations.
Multiplying \eqref{e:refdiff} with  $\langle \cdot,\psi \rangle_{H^{-1}}$ yields for $\rho\ge1$
\begin{align*}
\frac12\partial_t|\psi|_{H^{-1}}^2+|\psi|_{H^1}^2 
&\le |DF(u)[\psi]|_{L^1}|\psi|_{L^\infty}\\
&\le  C |\psi|_{L^\infty}\Big(|\psi|_{L^1}+|u|_{H^1}|\psi|_{H^1}\Big)\chi_\rho(|u|^2_{H^1})\\
&\quad+C |\psi|_{L^\infty}\Big(|u|_{L^1}+|u|_{H^1}^2\Big)\chi_\rho'(|u|^2_{H^1})|u|_{H^1}|\psi|_{H^1}\\
&\le C \rho |\psi|_{H^1}|\psi|_{L^\infty}.
\end{align*}
where we used H\"older, Sobolev embedding and the definition of the cut-off $\chi$.

Using Sobolev embedding of $L^\infty$ into $H^{\delta}$ for some 
$\delta>\frac12$ together with interpolation and Young inequality yields 
for some sufficiently large $p>1$ and some constant $c>0$
$$
\partial_t|\psi|_{H^{-1}}^2+|\psi|_{H^1}^2 \le \frac12|\psi|_{H^1}^2 +  c \rho^p |\psi|_{H^{-1}}^2.
$$
First, by  Gronwall  Lemma
$$
|\psi(t)|_{H^{-1}}^2\le |\psi(0)|_{H^{-1}}^2\mathrm{e}^{ct\rho^p},
$$ 
and then 
$$
\int_0^t|\psi|_{H^1}^2 dt
\le |\psi(0)|_{H^{-1}}^2 + c \rho^p \int_0^t  |\psi(s)|_{H^{-1}}^2\,ds
\le |\psi(0)|_{H^{-1}}^2\mathrm{e}^{ct\rho^p}.
$$
This together with \eqref{e:SFso} and the assumption on $\mathcal{Q}$ finishes the proof.
\end{proof}
\subsection{Some consequences}
It is well known that the strong Feller property implies that the laws
$P(t,x,\cdot)$ are mutually equivalent, for all $x$ and $t$. A less
obvious fact, which follows from Theorem 13 of Flandoli \& Romito
\cite{FlRo06c}, is that the same property holds between different
selection. In details, if $P^{(1)}(t,x,\cdot)$ and
$P^{(2)}(t,x,\cdot)$ are the Markov kernels associated to two
different selections, then $P^{(1)}(t,x,\cdot)$ and
$P^{(2)}(t,x,\cdot)$ are mutually equivalent for all $x$ and $t$.

Before enumerating all other properties following from strong
Feller, we need to show a technical result on the support of
the measures $P(t,x,\cdot)$. Following Flandoli \& Romito
\cite{FlRo06b}, we say that a Borel probability measure
$\mu$ is \emph{fully supported} on $H^1_{per}$ if $\mu[A]>0$
for every open set $A$ in $H^1_{per}$.
\begin{proposition}[Support theorem]\label{p:support}
Under Assumption \ref{a:noiseass}, let $(P_x)_{x\in L^2_{per}}$
be an a.\ s.\ Markov family. For every $x\in H^1_{per}$ and
$T>0$ the image measure of $P_x$ at time $T$ is fully supported
on $H^1_{per}$.
\end{proposition}
\begin{proof}
The proof is rather technical but straightforward, we only
give a sketch of it. To this purpose, we follow the same
steps of Flandoli \cite{Fla95} (see also Proposition 6.1
of \cite{FlRo06b}). It turns out that, since by Assumption
\ref{a:noiseass} the Wiener measure driving the equation
is fully supported on suitable spaces, we only have to
analyse the following control problem
\begin{equation}\label{e:control}
\dot h+h_{xxxx}=[-h_{xx}+(h_x^2)_{xx}]\chi_\rho+\dot w,
\qquad h(0)=x,
\end{equation}
where $w$ is the control. More precisely, we need to prove
the following two statements.
\begin{enumerate}
\item Given $T>0$, there is $\lambda\in(0,1)$ such that for $\rho>0$, $x\in H^1_{per}$,
$y\in H^4_{per}$ with $|x|_{H^1}\le\lambda\rho$ and $|y|_{H^1}\le\lambda\rho$, there
are $w\in\text{Lip}([0,T];H^1_{per})$ and $h\in C([0,T];H^1_{per})$ that solve \eqref{e:control}
with $h(T)=y$ and $\tau_\rho(w)>T$, where $\tau_\rho$ is defined as in \eqref{e:blowuptime}.
\item Let $w_n\to w$ in $W^{s,p}([0,T];D(A^\beta))$, with $s\in(\frac38, \frac12)$, $p>1$
such that $sp>1$ and $\beta\in(\frac14-s,-\frac18)$. Let $h_n$, $h$ be the solutions to
\eqref{e:control} corresponding to $w_n$, $w$ and let $\tau_n=\tau_\rho(w_n)$ and
$\tau=\tau_\rho(w)$. If $\tau>T$, then $\tau_n>T$ for sufficiently large $n$ and $h_n\to h$
in $C([0,T];H^1_{per})$.
\end{enumerate}
For the first claim, one uses \eqref{e:regmild} with $w=0$ to get a time $T_*<T$ such that
$h(T_*)\in H^4$ and $|h(T_*)|_{H^1}\le\rho$ (here we choose $\lambda$, using the estimates
on the semigroup $\textrm{e}^{tA}$). Then $h$ is given in $[T_*,T]$ by linear interpolation
from $h(T_*)$ to $y$ and $w$ in such a way that \eqref{e:control} is satisfied.

For the second claim, $s$, $p$ and $\beta$ are chosen so that the Wiener measure corresponding
to the random perturbation gives probability $1$ to $W^{s,p}([0,T];D(A^\beta))$ and the
convergence of $w_n$ implies that $z_n\to z$ in $C([0,T];H^1_{per})$, where $z_n$, $z$ are
the solutions to $\dot z=-z_{xxxx}+\dot w$ corresponding to $w_n$ and $w$ (this also gives
a common bound to $\tau_n$ and $\tau$, as in Lemma \ref{lem:SF-stop}). From this, it is easy
to see, by the mild formulation \eqref{e:regmild}, that $h_n\to h$.
\end{proof}
\begin{proposition}[Local regularity]
Let $(P_x)_{x\in L^2_{per}}$ be an a.\ s.\ Markov family and assume Assumption \ref{a:noiseass}.
Then for each $x\in H^1_{per}$ and all times $t>0$,
$$
P_x[\text{there is $\varepsilon>0$ such that }\xi\in C((t-\varepsilon,t+\varepsilon);H^1_{per})]=1.
$$
Moreover, for each $x\in H^1_{per}$, the set $T_{P_x}$ of property {\rm\textbf{\footnotesize\textsf{[E3]}}}
is empty, that is the energy inequality holds for all times.
\end{proposition}
\begin{proof}
Let $(\mathcal{P}_t)_{t\ge0}$ be the transition semigroup defined
by the given Markov family and set $\widetilde{\nu}=\int_0^1(\mathcal{P}_s^*\delta_0)\,ds$.
Set moreover $\widetilde{\Omega}_{a,b}=\{\xi\in C((a,b);H^1)\}$ and
$\widetilde{\Omega}_t=\bigcup\widetilde{\Omega}_{t-\varepsilon,t+\varepsilon}$.
We first observe that by \eqref{e:logtightforh},
$$
\widetilde{P}[\,|\xi_t|^2_{H^1}\ge\rho]
=\int_{t}^{t+1}P_0[\,|\xi_s|^2_{H^1}\ge\rho]\,ds
\le\frac{C}{\log(1+\rho)},
$$
where in particular the constant $C$ depends on $t$ (but it is increasing in $t$). Now, by
the Markov property, for all $\rho>0$,
$$
\widetilde{P}[\widetilde{\Omega}_{t-\varepsilon,t+\varepsilon}]
=\int P_y[\widetilde{\Omega}_{\varepsilon,3\varepsilon}]\,\pi_{t-2\varepsilon}\widetilde{P}(dy)
\ge(\inf_{|y|_{H^1}\le\rho}P_y[\widetilde{\Omega}_{\varepsilon,3\varepsilon}])(1-\frac{C}{\log(1+\rho)}),
$$
where $\pi_{s}\widetilde{P}$ is the marginal of $\widetilde{P}$ at time $s$. By Lemma \ref{lem:SF-stop}
we know that $\inf_{|y|_{H^1}\le\rho}P_y[\widetilde{\Omega}_{\varepsilon,3\varepsilon}]\uparrow1$
as $\varepsilon\to0$ and in conclusion $\widetilde{P}[\widetilde{\Omega}_t]=1$.

By disintegration, $P_x[\widetilde{\Omega}_t]=1$ for $\widetilde{\nu}$-a.\ e.\ $x$, hence for
a dense set of $H^1_{per}$ by Proposition \ref{p:support} and in conclusion for all $x\in H^1_{per}$
by the strong Feller property.
\end{proof}
The previous proposition and Theorem 6.7 of \cite{FlRo06b} (suitably adapted
to this framework) improve our knowledge on the Markov property as follows.
\begin{cor}
Under Assumption \ref{a:noiseass}, if $(P_x)_{x\in L^2_{per}}$ is an a.\ s.\ Markov
family of solutions to \eqref{e:dirkeq}, then $(P_x)_{x\in H^1_{per}}$ is a Markov
process. Namely
$$
\mathbb{E}^{P_x}[\varphi(\xi_t)|\mathcal{B}_s]=\mathbb{E}^{P_{\xi_s}}[\varphi(\xi_{t-s})],
\qquad P_x-a.s.,
$$
for all $x\in H^1_{per}$, $\varphi\in C_b(L^2_{per})$ and $0\le s\le t$. 
\end{cor}
\section{Existence and uniqueness of invariant measures}\label{s:ergodic}
Existence of an invariant measure for \eqref{e:dirkeq} is straightforward for
trace-class noise, as one can rely on It\^o formula applied to the energy
balance given by $|h(t)|_{L^2}^2$. The standard approximation is then tight,
since we can control $\mathbb{E}\bigl[\int_0^T|h_{xx}|_{L^2}^2\,dt\bigr]$.

In this section we prove existence of an invariant measure for more general
noise (such as space time white noise) under the assumption (which will be
valid for the whole section) that the equation has no linear instability,
namely
\begin{equation}\label{e:dirkeqstable}
\dot h = -h_{xxxx} + (h_x)^2_{xx}+\eta.
\end{equation}
In order to take the linear instability into account, gauge functions have to
be used, as in Bl\"omker \& Hairer \cite{BlHa04} or Collet et al. \cite{Eckmann},
Temam \cite{Tem2}, but up to now this is quite technical and only applicable
to Dirichlet or Neumann boundary conditions. For periodic boundary conditions
this question is still open.
\begin{theorem}\label{t:exIM}
Let $(P_x)_{x\in L^2_{per}}$ be any a.\ s.\ Markov family of energy martingale
solutions to \eqref{e:dirkeqstable}. Then there exists an invariant measure for
the transition semigroup associated to $(P_x)_{x\in L^2_{per}}$ with support
contained in $H^\gamma_{per}$, for some $\gamma\in(\frac54,\frac32)$.
\end{theorem}
\begin{remark}
Note that the upper bound $\gamma<\frac32$ is stated only for convenience.
The crucial restriction is  $\gamma>\frac54$, as in the proof of this
theorem we shall need  that $\widetilde{Z}_{\alpha,\cdot}\in W^{1,4}$,
which is implied by $\widetilde{Z}_{\alpha,\cdot}\in H^\gamma_{per}$,
where $\widetilde{Z}_{\alpha,\cdot}$ is the process defined in \eqref{e:ztilde}.
\end{remark}
By the results of the previous section we can immediately conclude that the
invariant measure is unique (via the strong Feller property and Doob's
theorem) and that it is fully supported on $H^1_{per}$ (by means of Proposition
\ref{p:support}).
\begin{cor}
Under Assumption \ref{a:noiseass}, the invariant measure provided by Theorem
\ref{t:exIM} above is unique and fully supported on $H^1_{per}$.
\end{cor}
So far we know that each Markov solution has its own unique invariant
measure. In principle, these invariant measures come from different
transition semigroups and do not need to be equal, even though they
have something in common. For example, we know from \cite[Theorem 13]{FlRo06c}
that they are mutually equivalent. At this stage, the problem of uniqueness
of the invariant measure over all selection is open, as well as the
well posedness of the martingale problem.
\subsection{The proof of Theorem \ref{t:exIM}}\label{ss:exIM}
Consider the family of measures of the Krylov-Bogoliubov method starting from
the initial condition $0$,
$$
\mu_T=\frac1T\int_0^T P_0[\xi(s)\in\cdot]\,ds.
$$
It is sufficient to prove compactness of $(\mu_{T})_{T\in \mathbf{N}}$ in $H^\gamma$.
Thus we need that for all $\varepsilon>0$ there is $R>0$ such that
\begin{equation}\label{e:aim}
\mu_{T}[\,|\cdot|_{H^\gamma}>2R]<\varepsilon,
\qquad\text{for all }T\in \mathbf{N}.
\end{equation}
First we consider $\widetilde{V}=\xi-\widetilde{Z}_{\alpha,z_0}$ 
for any initial condition $z_0\in H^\gamma$. As in Remark
\ref{r:eqweak}, $\widetilde{V}$ satisfies for some $\alpha>0$ 
$$
\dot{\widetilde{V}}+\widetilde{V}_{xxxx}
=\bigl[(\widetilde{V}_{x}+\widetilde{Z}_{x})^2\bigr]_{xx}+\alpha\widetilde{Z},
\qquad \widetilde{V}(0)=-z_0
$$
and $\widetilde{Z}=\widetilde{Z}_{\alpha,z_0}$ is a solution of
$$
d\widetilde{Z}+(\alpha\widetilde{Z}+\widetilde{Z}_{xxxx})\,dt=dW,
\qquad \widetilde{Z}(0)=z_0.
$$
Now we can bound
\begin{align*}
P_0[|\xi(s)|_{H^{\gamma}}>2R]
&\le P_0[\,|\widetilde{V}(s)|_{H^\gamma}+|\widetilde{Z}(s)|_{H^\gamma}>2R]\\
&\le P_0[\,|\widetilde{V}(s)|_{H^\gamma}>R] + P_0[|\widetilde{Z}(s)|_{H^\gamma}>R]\\
&\le P_0[\,|\widetilde{V}_{xx}(s)|>R] + P_0[|\widetilde{Z}(s)|_{H^\gamma}>R].
\end{align*}

Let $\varphi:[0,\infty)\rightarrow\mathbb{R}$ be a function,
 which we will determine at the end of the proof, such that
$\varphi$ is increasing, concave, with $\varphi(r)\uparrow\infty$
as $r\uparrow\infty$, and for every $x$, $y\ge0$,
\begin{equation}\label{e:condphi}
\varphi( x+y)
\le C+\varphi(x)+\log(y+1)
\le C+\varphi(x)+y.
\end{equation}
Note that we are not able to bound moments or log-moments of $\widetilde{V}$
uniformly in time. All we can show is that the $\varphi$-moment is bounded
uniformly in time\footnote{Bl\"omker \& Hairer \cite{BlHa04} give a different
proof of existence of an invariant measure, which relies on Galerkin approximations.
Here we consider any solution to the equation, which in principle could not be
a limit of such approximations, if the solutions are not unique.}. Consider
\begin{align}\label{e:exsys1}
\frac1T\int_0^T P_0[\,|\widetilde{V}_{xx}(s)|_{L^2}>R]  
&=   \frac1T\int_0^T P_0[\varphi(|\widetilde{V}_{xx}(s)|^2_{L^2})>\varphi(R^2)]\,ds\notag\\
&\le \frac{1}{\varphi(R^2)}\mathbb{E}^{P_0}\Bigl[\frac1T\int_0^T\varphi(|\widetilde{V}_{xx}(s)|^2_{L^2})\,ds\Bigr].
\end{align}
From the energy  inequality we know that for all $t$ and almost every $t_0\in[0,t]$,
\begin{align*}
|\widetilde{V}(t)|_{L^2}^2+\int_{t_0}^t|\widetilde{V}_{xx}(s)|_{L^2}^2\,ds 
&\le    |\widetilde{V}(t_0)|_{L^2}^2+C\int_{t_0}^t|\widetilde{Z}_x(s)|_{L^4}^{\frac{16}{3}}|\widetilde{V}(s)|_{L^2}^2\,ds\\
&\quad +C\int_{t_0}^t(|\widetilde{Z}_x(s)|_{L^4}^4+\alpha^2|\widetilde{Z}(s)|_{L^2}^2)\,ds
\end{align*}
Let us fix some notation:
$$
a(t)=C|\widetilde{Z}_x(t)|_{L^4}^{\frac{16}{3}},
\qquad
b(t)=C(|\widetilde{Z}_x(t)|_{L^4}^4+\alpha^2|\widetilde{Z}(t)|_{L^2}^2),
$$
where all moments of $a$ and $b$ are bounded by some constant and the initial
condition $\widetilde{Z}(0)$. Thus for all $t>0$ and almost all $t_0\in[0,t]$,
\begin{equation}\label{e:exener1}
|\widetilde{V}(t)|_{L^2}^2+\int_{t_0}^t|\widetilde{V}_{xx}(s)|_{L^2}^2\,ds
\le|\widetilde{V}(t_0)|_{L^2}^2+\int_{t_0}^t(a(s)|\widetilde{V}(s)|_{L^2}^2+b(s))\,ds.
\end{equation}
Using Poincar\'e inequality (with constant $\lambda$) it follows that
\begin{equation} \label{e:exener2}
|\widetilde{V}(t)|_{L^2}^2+\int_{t_0}^t(\lambda-a(s))|\widetilde{V}(s)|_{L^2}^2\,ds
\le|\widetilde{V}(t_0)|_{L^2}^2+\int_{t_0}^t b(s)\,ds.
\end{equation}
We now compare $\widetilde{V}$  with some simpler one-dimensional equation.
Let $u(t)$ be the solution of 
$$
u(t)+\int_{0}^{t}(\lambda-a(s))u(s)\,ds
=|z_0|_{L^2}^2+\int_{0}^{t}b(s)\,ds,
$$
namely,
$$
u^{\prime}(t)  +(\lambda-a(t))u(t)  =b(t),
\qquad u(0)=|z_0|_{L^2}^2.
$$
This is exactly the situation of the modified Gronwall Lemma \ref{p:gronwall},
hence we derive $|\widetilde{V}(t)|_{L^{2}}^{2}\leq u(t)$ and thus also
\begin{align}\label{e:ex1}
\int_0^t|\widetilde{V}_{xx}(s)|_{L^2}^2\,ds
&\le u(0)+\int_0^t(a(s)u(s)+b(s))\,ds\notag\\
&\le u(0)+\int_0^t(a^2(s)+u_*(s)+b(s))\,ds,
\end{align}
where $u_*(t)$ is the solution to the one-dimensional equation
\begin{equation}\label{e:defu*}
u_*'(t)+(\lambda-a_*(t))u_*(t)=b_*(t)
\end{equation}
with
$$
a_*(t)=2a(t),\qquad
b_*(t)=Cb^2(t),\qquad
u_*(0)=u^2(0).
$$
Note that, as
$$
(u^2(t))'+2(\lambda-a(t))u^2(t)
=  2u(t)b(t)
\le\lambda u^2(t)+Cb^2(t),
$$
with a constant depending on $\lambda$, we have by
a comparison principle for ODEs that $u^2(t)\le u_*(t)$.

Let us consider for notational simplicity only the case of integer
$T$. From \eqref{e:exsys1}
\begin{equation*}
\frac1T\int_0^T P_{0}[\,|\widetilde{V}_{xx}(s)|^2_{L^2}>R]
\le\frac{1}{\varphi(R^2)}\frac1T\sum_{k=0}^{T-1}\mathbb{E}^{P_0}\Bigl[\varphi\bigl(\int_{k}^{k+1}|\widetilde{V}_{xx}(s)|^2_{L^2}\,ds\bigr)\Bigr].
\end{equation*}
Thus we only need to bound these moments independently of $k$. 
The splitting in discrete time steps is necessary, in order to avoid 
suprema over $[0,T]$, which usually give $T\log(T)$ terms.

From \eqref{e:exener1},
\begin{align*}
\int_{k}^{k+1} |\widetilde{V}_{xx}(s)|_{L^2}^2\,ds
&\le |\widetilde{V}(k)|_{L^2}^2 + \int_{k}^{k+1}(a(s)|\widetilde{V}(s)|_{L^2}^2 + b(s))\,ds\\
&\le \sup_{s\in[k,k+1]}|\widetilde{V}(s)|_{L^2}^4 + \int_{k}^{k+1}(a^2(s)+b(s))\,ds +1.
\end{align*}
We use the fact that the stochastic convolution is bounded in
expectation by a constant plus the initial condition, i.e.\
$\mathbb{E}^{P_0}|\widetilde{Z}(t)|^p\le C(1+|z_0|^p)$ in $L^2$,
$H^\gamma$, and $W^{1,4}$. We derive
\begin{align*}
\mathbb{E}^{P_0}\bigl[\varphi\bigl(\int_{k}^{k+1}|\widetilde{V}_{xx}(s)|_{L^2}^2\,ds\bigr)\bigr]
&\le C+\mathbb{E}^{P_0}\Bigl[\varphi\bigl(\sup_{s\in[k,k+1]}|\widetilde{V}(s)|_{L^2}^4\bigr)\Bigr]\\
&\quad+\mathbb{E}^{P_0}\Bigl[\int_{k}^{k+1}(a^2(s)+b(s)+1)\,ds\Bigr]\\
&\le C+|z_0|_{W^{1,4}}^{11}+\mathbb{E}^{P_0}\bigr[\varphi\bigl(\sup_{s\in[k,k+1]}|\widetilde{V}(s)|_{L^2}^4\bigr)\bigr].
\end{align*}
Thus
\begin{multline*}
\frac1T\int_0^T P_0[\,|\widetilde{V}_{xx}(s)|^2_{L^2}>R]\le\\
\le \frac1{\varphi(R^2)}\Bigl(C+|z_0|_{W^{1,4}}^{11}+\sup_{k\le T-1}\mathbb{E}^{P_0}\bigl[\varphi(\sup_{s\in[k,k+1]}|\widetilde{V}(s)|_{L^2}^4)\bigr]\Bigr).
\end{multline*}
As $|\widetilde{V}(t)|_{L^2}^4\le u^2(t)\le u_*(t)$, it is sufficient to show
that there are a function $\varphi$ with all the above specified properties
and a constant $C>0$ such that
$$
\mathbb{E}^{P_0}\bigl[\varphi(\sup_{t\in[k,k+1]}u_*(t))\bigr]\le C
$$
independently of $k$. Recall the choice $\xi(0)=0$, and thus $\widetilde{V}(0)=-z_0$
in the Krylov-Bogoliubov scheme. From \eqref{e:defu*}, $u_*(t)$ is given by
\begin{equation}\label{e:bouu*}
u_*(t)  
=\int_0^t\mathrm{e}^{\int_s^t(-\lambda+a_*(r))\,dr}b_*(s)\,ds
 +\mathrm{e}^{\int_s^t(-\lambda+a_*(r))\,dr}|z_0|_{L^2}^4.
\end{equation}
Recall the special shape of $a_*$ and $b_*$. By renaming constants,
$$
a_*(t) =C_*|\widetilde{Z}_x(t)|_{L^4}^{\frac{16}{3}},
\qquad
b_*(t)\le C_*|\widetilde{Z}_x(t)|_{L^4}^8+C_*\alpha^2|\widetilde{Z}(t)|_{L^2}^4.
$$
Set moreover
$$
\theta(t)=C_*|\widetilde{Z}_x(t)|_{L^4}^{\frac{16}{3}}
         +C_*|\widetilde{Z}_x(t)|_{L^4}^8
         +C_*|\widetilde{Z}(t)|_{L^2}^4.
$$
As $\mathbb{E}^{P_0}[\theta(t)]\to 0$ for $\alpha\to0$,
we choose $\alpha$ sufficiently large such that
$$
\mathbb{E}^{P_0}[\theta(t)]\le\frac{\lambda}{4}.
$$
From \eqref{e:bouu*}
\begin{align*}
u_*(t)
&  \le  (1+\alpha^2)\int_0^t\mathrm{e}^{\int_s^t(-\lambda+\theta(r))\,dr}\theta(s)\,ds
       +\mathrm{e}^{\int_0^t(-\lambda+a_*(r))\,dr}|z_0|_{L^2}^4\\
&    =  \lambda(1+\alpha^2)\int_0^t\mathrm{e}^{\int_s^t(-\lambda+\theta(r))\,dr}\,ds
       +\mathrm{e}^{\int_0^t(-\lambda+a_*(r))\,dr}|z_0|_{L^2}^4\\
&\quad +(1+\alpha^2)\int_0^t\mathrm{e}^{\int_s^t(-\lambda+\theta(r))\,dr}(-\lambda+\theta(s))\,ds\\
&  \le  \lambda(1+\alpha^2)\int_0^t\mathrm{e}^{\int_s^t(-\lambda+\theta(r))\,dr}\,ds
       +(1+\alpha^2+|z_0|_{L^2}^4)\mathrm{e}^{\int_0^t(-\lambda+\theta(r))\,dr}.
\end{align*}
Denote by $u_{**}(t)$ the function
$$
u_{**}(t)=\int_0^t\mathrm{e}^{\int_s^t(-\lambda+\theta(r))\,dr}\,ds
$$
which is a solution of
$$
u_{**}'(t)+(\lambda-\theta(t))u_{**}(t)=1,
\qquad u_{**}(0)=0.
$$
Then
\begin{align*}
\varphi(\sup_{t\in[k,k+1]}u_*(t))
&\le  C+\varphi\bigl(\lambda(1+\alpha^2)\sup_{t\in[k,k+1]}u_{**}(t)\bigr)\\
&\quad +\log\bigl((1+\alpha^2+|z_0|_{L^2}^4)\mathrm{e}^{\int_0^t(-\lambda+\theta(r))\,dr}\bigr).
\end{align*}
Thus bounding the stochastic convolution
\begin{multline*}
\frac1T\int_0^TP_0[\,|\widetilde{V}_{xx}(s)|^2_{L^2}>R]\,ds\le\\
\le\frac{1}{\varphi(R^2)}\Bigl(C+|z_0|_{W^{1,4}}^{11}
     +\sup_{k\in\mathbb{N}}\mathbb{E}^{P_0}\bigl[\varphi\bigl(\lambda(1+\alpha^2)\sup_{t\in[k,k+1]}u_{**}(t)\bigr)\bigr]\Bigr).
\end{multline*}
Let us now turn to bound $u_{**}$,
\begin{align*}
u_{**}(t)    
&\le \sup_{s\in[0,t]}\mathrm{e}^{\int_s^t(-\frac\lambda2+\theta(r))\,dr}\int_0^t\mathrm{e}^{-\frac\lambda2(t-s)}\,ds\\
&\le \frac{2}{\lambda}\exp\bigl[\sup_{s\in[0,t]}\int_s^t(-\frac{\lambda}{2}+\theta(r))\,dr\bigr],
\end{align*}
hence we need to bound
$$
\sup_{k\ge0}\mathbb{E}^{P_0}\Bigl[\varphi\bigl(2(1+\alpha^2)\exp(\sup_{t\in[k,k+1]}\sup_{s\in[0,t]}\int_s^t(-\frac{\lambda}{2}+\theta(r))\,dr)\bigr)\Bigr]<\infty.
$$
But, as we have
$$
\int_s^t(-\frac{\lambda}{2}+\theta(r))\,dr
\le\int_s^{\lceil t\rceil}(-\frac{\lambda}{2}+\theta(r))\,dr+\frac{\lambda}{2},
$$
we derive
$$
\sup_{t\in[k,k+1]}\sup_{s\in[0,t]}\int_s^t(-\frac{\lambda}{2}+\theta(r))\,dr
\le\frac{\lambda}{2}+\sup_{s\in[0,k+1]}\int_s^{k+1}(-\frac{\lambda}{2}+\theta(r))\,dr
$$
and thus we finally obtain  for $T\in\mathbf{N}$,
\begin{multline*}
\frac1T\int_0^T P_0[\,|\xi(s)|^2_{H^\gamma}>R]\,ds\le\\
\le\frac{1}{\varphi(R^2)}\Bigl(C+|z_0|_{W^{1,4}}^{11}
       +\sup_{k\in\mathbb{N}}\mathbb{E}^{P_0}\bigl[\varphi\bigl(C_{\alpha,\lambda}\exp(\sup_{s\in[0,k+1]}\int_s^{k+1}\!\!(-\frac{\lambda}{2}+\theta(r))\,dr)\bigr)\bigr]\Bigr).
\end{multline*}
Now we can use Lemma \ref{lem:OU-trick} to replace the 
OU-process in $\theta$ by the stationary process, thus
obtaining a process $\widetilde{\theta}$. Furthermore,
$z_0$ is replaced by $\widetilde{Z}(0)$.
Note that $\{\widetilde{\theta}(t)\}_{t\in\mathbf{R}}$ 
is now no longer defined on the same probability space
as $\theta$. Thus the expectation also changes.
Due to stationarity we have
\begin{align*}
\sup_{s\in[ 0,k+1]}\int_s^{k+1}(-\frac{\lambda}{2}+\tilde\theta(r))\,dr
&\overset{(L)}{=} \sup_{s\in[-(k+1),0]}\int_s^0(-\frac{\lambda}{2}+\widetilde{\theta}(r))\,dr\\
&\le              \sup_{s\in(-\infty,0]}\int_s^0(-\frac{\lambda}{2}+\widetilde{\theta}(r))\,dr\\
&\overset{(L)}{=} \sup_{t\in[0,\infty)}\int_0^t(-\frac{\lambda}{2}+\widetilde{\theta}(r))\,dr.
\end{align*}
Therefore, if we define the random variable
$$
\widetilde{X}=\sup_{t\in[0,\infty)}\int_0^t(-\frac{\lambda}{2}+\widetilde{\theta}(r))\,dr,
$$
we only have to prove that there exists a function $\varphi$ as above such that
$$
\widetilde{\mathbb{E}}[\varphi(2(1+\alpha^2)\mathrm{e}^{\widetilde{X}})]<\infty.
$$
Since $\widetilde{X}$ is finite with probability one by the ergodic theorem, such a
$\varphi$ exists by Lemma \ref{p:moment}. The proof of Theorem \ref{t:exIM} is
complete.
\begin{remark}
In the previous proof, we were only able to bound some moment of $V_{xx}$,
but using the trick of Debussche \& Da Prato \cite{deb-trick}, where
$\alpha$ is allowed to be random, it is possible to bound arbitrary polynomial
moments on bounded time intervals. 
\end{remark}
\begin{lemma}\label{lem:OU-trick}
Let $\delta>0$ and let $\phi$ be a positive map defined on the probability space $\Omega$.
If for all $z_0$
$$
\delta\le\mathbb{E}^{P_x}[\phi(\widetilde{Z}_{\alpha,z_0})]
\qquad \text{for $P_x$-almost every $\xi\in\Omega$},
$$
where $\widetilde{Z}_{\alpha,z_0}$ is the Ornstein-Uhlenbeck
process starting in $z_0$, as defined in \eqref{e:ztilde},
then
$$
\delta\le\int_{H^1}\phi(z)\,\mu_{OU}^*(dz),
$$
where $\mu_{OU}^*$ is the law of the stationary Ornstein-Uhlenbeck process.
\end{lemma}
The lemma is easily proved by averaging both sides with respect to
$z_0$ with the stationary Ornstein-Uhlenbeck process and using Tonelli
theorem.
\section{A priori estimates}\label{sec:WMS}
In this section we state all regularity results on processes
$Z$ and $V$. The first part contains the results on $Z$ under an
arbitrary weak martingale solution (from Definition \ref{d:weakms}).
Similarly, the second part contains the results on $V$ under an
arbitrary energy martingale solution (from Definition \ref{d:energyms}).
\subsection{Weak martingale solution}
Here we will present some lemmas on the regularity of $Z$ without
using equivalent versions, since our approach forces us to keep
the canonical process.
\begin{lemma}
\label{lem:ZregL4}
Given a weak martingale solution $P$, then for every $T>0$,
$$
\mathbb{E}^P\int_0^T|Z(t)|^4_{W^{1,4}}\,dt<\infty.
$$
\end{lemma}
\begin{proof}
It is enough to verify that $(Z_x)^2 \in L^2(\Omega\times(0,T),L^2_{per})$.
From the definition, we can write  $Z(t)$ as a complex Fourier series,
such that
$$
Z_x=\sum_{k\not=0} I_k \mathrm{e}^{ikx},
$$
where $I_k$ is a time dependent Gaussian real valued random variable
with $\mathbb{E}^P I_k^2 \le C |k|^{-2}$. Thus,
$\mathbb{E}^P I_k^4 \le C |k|^{-4}$, too. Now,
$$
(Z_x)^2=
\sum_{n\in\mathbf{Z}}\sum_{k\not=0,n} I_k I_{n-k}  \mathrm{e}^{inx},
$$
We derive
\begin{align*}
\mathbb{E}^P |(Z_x)^2|_{L^2}^2
&=   \sum_{n\in\mathbf{Z}}\mathbb{E}^P\Big(\sum_{k\not=0,n} I_k I_{n-k}\Big)^2\\
&\le \sum_{n\in\mathbf{Z}}\sum_{k\not=0,n}\sum_{l\not=0,n}\mathbb{E}^P[|I_k||I_{n-k}||I_l| |I_{n-l}|]\\
&\le \sum_{n\in\mathbf{Z}}\Big(\sum_{k\not=0,n}\frac1{|k||n-k|}\Big)^2,
\end{align*}
where we used H\"older's inequality in the last step. It is an
elementary exercise to check that the series in the last equation
converges. Thus integration in time yields the result.
\end{proof}
\begin{lemma}\label{lem:ZregW14}
Let $P$ be a weak martingale solution. Then
for some $\lambda>0$ there is a constant $C$ such that
$$
\int_0^T\mathbb{E}^{P}\exp\{\lambda\|Z_x(t)\|^2_{L^4} \}\,dt \le CT
\qquad\text{for all $T>0$}.
$$
Thus, for some constant $C$ depending only on $q$, $p$ and $T$,
$$
\sup_{T\ge0}\frac1T\mathbb{E}^{P}\|Z\|^p_{L^p([0,T],W^{1,4})}<\infty
\qquad\text{and}\qquad
\sup_{T\ge0}\frac1T\mathbb{E}^{P}\|Z\|^q_{L^p([0,T],W^{1,4})}\le C.
$$
\end{lemma}
\begin{proof}
Using Lemma \ref{lem:ZregL4} we know that $\mathbb{E}^{P}\|Z_x(t)\|^4_{L^4}\le C$
for all $t\ge0$. As $Z_x(t)$ is a Gaussian random variable in $L^4$,
Ferniques theorem (see Da Prato \& Zabczyk \cite{DPZa}) implies that
$$
\sup_{t\ge0}\mathbb{E}^{P}\exp\{\lambda\|Z_x(t)\|^2_{L^4} \}<\infty,
$$
for some $\lambda>0$. Thus
$$
\mathbb{E}^{P}\|Z\|^p_{L^p([0,T],W^{1,4})}
\le C\int_0^T \mathbb{E}^{P}\exp\{\lambda\|Z_x(t)\|^2_{L^4} \}\,dt
\le C T,
$$
where the constant does not depend on $T$.
The last claim follows from H\"older inequality.
\end{proof}
The following lemma on the $L^\infty([0,\infty),L^2_{per})$-regularity
is necessary to transfer weak continuity in $L^2$ from $V$ to $Z$.
Note, again, that we cannot prove continuity of $Z$, as we are not
using continuous versions of the canonical process $Z$.
\begin{lemma}\label{lem:ZregLinfty}
Let $P$ be a weak martingale solution. Then for $0\le s<\frac32$, 
$p>1$ and $T>0$
$$
Z\in L^p(\Omega,L^\infty([0,T],H^s_{per}))
$$
and thus
$$
P[Z\in L^\infty([0,\infty),H^s_{per})]=1.
$$
Due to $Z\in\Omega$, we thus have $Z$ is $P$-a.s.\ weakly continuous
with values in $H^s_{per}$.
\end{lemma}
\begin{proof}
Using the factorisation method (see Da Prato \& Zabczyk \cite[Chapter 5]{DPZa}),
$$
Z(t)=C_\alpha \int_0^t \mathrm{e}^{(t-\tau)A}(t-\tau)^{\alpha-1}Y(\tau)\,d\tau,
$$
with
$$
Y(\tau)=\int_0^\tau\mathrm{e}^{(\tau-s)A}(\tau-s)^{-\alpha}\,dW(s).
$$
We fix $T>0$, $\alpha\in(0,\frac{3-2s}8)$ and $m>\frac1\alpha>\frac83$,
and let the constants depend on them. Now using H\"older's inequality,
$$
\sup_{t\in[0,T]}|Z(t)|_{H^s}
\le C \sup_{t\in[0,T]}\int_0^t (t-\tau)^{\alpha-1}|Y(\tau)|_{H^s}\,d\tau
\le C \Big(\int_0^T |Y(\tau)|_{H^s}^m\,d\tau\Big)^{\frac1m}.
$$
Thus using that $Y$ is Gaussian,
$$
\mathbb{E}^P \sup_{t\in[0,T]}|Z(t)|_{H^s}^m
\le  C \int_0^T (\mathbb{E}^P|Y(\tau)|^2_{H^s})^{m/2}\,d\tau
\le C \Big(\sum_{k=1}^\infty k^{2s}\alpha_k^2|\lambda_k|^{2\alpha-1} \Big)^{m/2}.
$$
The last series converges, as $\alpha_k^2\le C$ and $\lambda_k\sim -k^4$.
Taking $T\in\mathbf{N}$ concludes the proof.
\end{proof}
\subsection{Energy martingale solution}
This part is devoted to the proof of the \emph{tightness} property
for sequences of energy martingale solutions, essentially by means
of bounds on the process $V$.
\begin{lemma}\label{l:tighttransfer}
Let $(P_n)_{n\in\mathbf{N}}$ be a family of energy Markov solutions. Then
the sequence of laws of $V$ under $P_n$ is tight in $L^2(0,T,H^1_{per})$,
if and only if $(P_n)_{n\in\mathbf{N}}$ is tight in $L^2(0,T,H^1_{per})$.

The same result is true for any space in which $Z$ is defined,
for example $C(0,T,H^{-4}_{per})$.
\end{lemma}
\begin{proof}
We prove only one direction, the other one is the same.
As $\mathrm{Law}_{Z,n}=P_n[Z\in\cdot]$ is by Definition \ref{d:weakms} and Lemma
\ref{lem:ZregLinfty} the law of the stochastic convolution in
$L^2([0,T],H^1_{per})$ and thus independent of $n$.
Hence, the family of measures  $(\mathrm{Law}_{Z,n})_{n\in\mathbf{N}}$
is tight in $L^2([0,T],H^1_{per})$.
Thus there is a compact subset $K_{\varepsilon,1}\subset L^2([0,T],H^1_{per})$
with $P_n[Z\in K_{\varepsilon,1}]>1-\varepsilon$.
Furthermore, by the tightness of $P_n[V\in\cdot]$, there is a compact set
$K_{\varepsilon,2}\subset L^2([0,T],H^1_{per})$ such that
$P_n[V\in K_{\varepsilon,2}]>1-\varepsilon$.

Define now the compact subset
$$
K_{\varepsilon,3}
=K_{\varepsilon,1}+K_{\varepsilon,2}
=\{u=u_1+u_2|\ u_i\in K_{\varepsilon,i}\},
$$
then by $\xi=V+Z$ we have 
$$
P_n\left[K_{\varepsilon,3}\right]
\ge P_n\left[Z\in K_{\varepsilon,1},\ V\in K_{\varepsilon,2}\right]
\ge 1-2\varepsilon,
$$
which concludes the proof.
\end{proof}
\begin{lemma}\label{l:tderV}
Let $P$ be an energy martingale solution. Then for all $T>0$
$$
\|\partial_t V\|_{L^2([0,T],H^{-3}_{per})}
\le C \|V\|_{L^2([0,T],H^2)}(1+\|V\|_{L^\infty([0,T],L^2)}) + C \|Z\|_{L^4([0,T],H^1)},
$$
$P$-almost surely, with constants independent of $P$.
\end{lemma}
\begin{proof}
From Remark \ref{r:eqweak}, we know that for $\varphi\in H^3_{per}$
with $|\varphi|_{H^3}=1$ we have
$$
\partial_t \langle V,\varphi \rangle =
-\langle V_{xx}+V,\varphi_{xx}\rangle
-\langle (V_x+Z_x)^2,\varphi_{xx} \rangle
$$
Thus using the embedding of $L^1$ into $H^{-1}$ and
an interpolation inequality,
\begin{align*}
|\partial_t V|_{H^{-3}}
&\le |V_{xx}|_{L^2}+|V|_{L^2}+|(V_x+Z_x)^2|_{L^1},\\
&\le |V|_{H^2}+C|V|_{L^2}|V|_{H^2}+2|Z|_{H^1}^2.
\end{align*}
Integrating the square in time yields the result.
\end{proof}
\begin{lemma}\label{l:vbound}
Let $(P_n)_{n\in\mathbf{N}}$ be a family of energy martingale solutions.
Define
$$
\Sigma(R)=\{u: \|u\|_{L^\infty([0,T],L^2)}<R
\qquad\text{and}\qquad
\|u\|_{L^2([0,T],H^2)}<R\}.
$$
Suppose that $P_n$ is started at a probability measure $\mu_n$ such that
$$
\int_{L^2_{per}}\bigl(\log(|x|_{L^2}+1)\bigr)^\kappa\,\mu_n(dx)\le K,
$$
for all $n\in\mathbf{N}$ and for some $\kappa>0$, then
$$
\sup_{n\in\mathbf{N}}P_n[V\in\Sigma(R)]\ge 1-\frac{C}{\log(1+R)^\kappa}.
$$
\end{lemma}
\begin{proof}
By property \textbf{\textsf{\footnotesize[E3]}}, we have that, $P_n$-almost surely,
\begin{align*}
&|V(t)|_{L^2}^2+\int_0^t|V_{xx}|^2_{L^2}\le\\
&\qquad\le |V(0)|_{L^2}^2 +\int_0^t(|V_x|_{L^2}^2+2|V_x|_{L^4}|Z_x|_{L^4}|V_{xx}|_{L^2}+|Z_x|^2_{L^4}|V_{xx}|_{L^2})\,dt\\
&\qquad\le |V(0)|_{L^2}^2 +\int_0^t \frac12|V_x|_{L^2}^2+C(1+|Z_x|_{L^4}^{16/3}) |V|^2_{L^2}+C|Z_x|^4_{L^4}\,dt,
\end{align*}
where we have used the Sobolev embedding of $H^1_{per}$ into $L^4_{per}$, interpolation, Young, 
and Poincar\'e inequalities. Now from Gronwall's inequality it follows that, for all $t\in[0,T]$,
\begin{equation}\label{e:vbound}
|V(t)|_{L^2}^2+\int_0^t|V_{xx}|^2_{L^2}
\le C \mathrm{e}^{C \|Z\|^{16/3}_{L^{16/3}([0,T],W^{1,4})}}\Big(|V(0)|_{L^2}^2+\|Z\|^4_{L^4([0,T],W^{1,4})}\Big)
\end{equation}
where the constants might depend on $T$. Applying $(\log(x+1))^\kappa$ and using
the inequality
$$
\log(x+y+1)^\kappa\le C(\log(x+1)^\kappa+\log(y+1)^\kappa)
\quad\text{for }x,y\ge0,
$$
leads to
$$
\mathbb{E}^{P_n}\Bigl[\sup_{t\in[0,T]}\log(1+|V(t)|_{L^2}^2)\Bigr]^\kappa\le C
$$
and
$$
\mathbb{E}^{P_n}\Bigl[\log(1+\int_0^t|V_{xx}(s)|^2_{L^2}\,ds)\Bigr]^\kappa \le C,
$$
where the constant is independent of $n$. Now Chebychev inequality yields the result.
\end{proof}
The main result of this section is:
\begin{theorem}\label{t:tightness}
Let $(P_n)_{n\in\mathbf{N}}$ be a family of energy martingale solutions
with each $P_n$ starting in $\mu_n$ and
$$
\int_{L^2_{per}}[\log(|x|_{L^2}+1)]^\kappa\,\mu_n(dx)\le K,
\qquad\text{for all }n\in\mathbf{N},
$$
for some $\kappa>0$ and $K>0$. Then $(P_n)_{n\in\mathbf{N}}$ is tight on
$\Omega\cap L^2([0,\infty),H^1_{per})$.

Furthermore, there is a constant depending only on $T>0$, $z_0\in H^1$, $K>0$, and $\kappa>0$,
such that
\begin{align}
&\mathbb{E}^{P_n}\Bigl[\log\bigl(1+\int_0^T|\xi_x(s)|^2_{L^2}\,ds\bigr)\Bigr]^\kappa
\le C,\label{e:logtightforh}\\
&\mathbb{E}^{P_n}\Bigl[\log\bigl(1+\int_0^T|V_{xx}(s)|^2_{L^2}\,ds\bigr)\Bigr]^\kappa
  +\mathbb{E}^{P_n}\bigl[\sup_{t\in[0,T]}\log(1+|V(t)|^2_{L^2})\bigr]^\kappa\le C.\label{e:logtightforv}
\end{align}
\end{theorem}
\begin{proof}
For the bounds on logarithmic moments of $V$ we use the bounds obtained at the end 
of the proof of the previous Lemma \ref{l:vbound}. Using the bounds on $Z$ from Lemma
\ref{lem:ZregLinfty} yields the bound on logarithmic moments of $\xi$.

For the tightness of the law of $V$ under $P_n$ we use Lemmas  \ref{l:vbound}
and \ref{l:tderV} for  the bound for $\partial_t V$,
together with the compact embeddings of
$H^1([0,T],H^{-3}_{per})$ into $C([0,T],H^{-4}_{per})$
and of $L^2([0,T],H^{2}_{per})\cap H^1([0,T],H^{-3}_{per})$ into
$L^2([0,T],H^{1}_{per})$ (see for example Temam \cite{Tem}).

For the tightness of $P_n$ we use Lemma \ref{l:tighttransfer} on
transfer of tightness in the spaces $L^2([0,\infty),H^1_{per})$
and $C([0,T],H^{-4}_{per})$.
\end{proof}
\section{Some useful technical tools}\label{s:tools}
\subsection{A suitable concave moment}
We aim to prove the following proposition.
\begin{proposition}
\label{p:moment}
Let $X$ be a random variable with values in $[0,\infty)$. Then there
is a concave and non-decreasing map $\phi:[0,\infty)\to[0,\infty)$
such that $\phi(x)\uparrow\infty$ and
$$
\mathbb{E}[\phi(X)]<\infty.
$$
Moreover, $\phi$ can be chosen in such a way that for some constant $C$,
$$
\phi(x+y)\le\phi(x)+Cy,
\qquad\text{for all }x,y\in[0,\infty).
$$
\end{proposition}
\begin{remark}
Notice that the last condition on $\phi$ given in the proposition
above can be replaced by
$$
\phi(x+y)\le\phi(x)+C\log(1+y)
$$
for some constant $C>0$ and for all $x$, $y\in[0,\infty)$. Indeed, 
let $\varphi$ be the map given by the proposition, then
$\phi(x)=\varphi(\log(1+x))$ has exactly the same properties of $\varphi$
and $\phi(x+y)=\varphi(\log(1+x+y))\le\phi(x)+C\log(1+y)$, since
$\log(1+x+y)\le\log(1+x)+\log(1+y)$.
\end{remark}
\begin{proof}
We first show that there is a non-decreasing continuous map $u:[0,\infty)\to[0,\infty)$
such that $u(0)=0$, $u(x)\uparrow\infty$ as $x\to\infty$ and $\mathbb{E}[u(X)]<\infty$.
Choose a sequence $(x_n)_{n\in\mathbf{N}}$ such that $x_0=0$, $x_n\uparrow\infty$
and $4^n\mathbb{P}[x_n\le X<x_{n+1}]\longrightarrow1$.
This can always been done, since $X$ is a.\ s.\ finite. Now, let $\widetilde{u}$
be the piece-wise constant function that on each interval $[x_n,x_{n+1})$ takes
the value $2^n$. We finally set
$u(t)=\frac{1}{t}\int_0^t[\widetilde{u}(t)-\inf_{s\ge0}\widetilde{u}(s)]\,ds$.

Next, we show how to construct a map $\phi$ as in the statement of the
proposition such that $\phi\le 1+u$. Define the sequence
$(y_n)_{n\in\mathbf{N}}$ as $y_0=0$ and $y_n=\max\{x\in[0,\infty):u(x)=n\}$,
for $n\ge1$. The sequence $(y_n)_{n\in\mathbf{N}}$ is increasing and
$y_n\uparrow\infty$. Define $\phi$ as $\phi(y_0)=0$, $\phi(y_1)=1$,
$$
\phi(y_n)=\min\bigl\{n,\phi(y_{n-2})+\frac{\phi(y_{n-1})-\phi(y_{n-2})}{y_{n-1}-y_{n-2}}(y_n-y_{n-2})\bigr\},
$$
and by linear interpolation for all other values of $x\in[0,\infty)$. In other
words, at each point $y_n$ the map is defined as either the continuation of the
line $y_{n-2}\longrightarrow y_{n-1}$ or $u(y_n)$, depending on which is
the smallest value. The construction is shown in the picture.
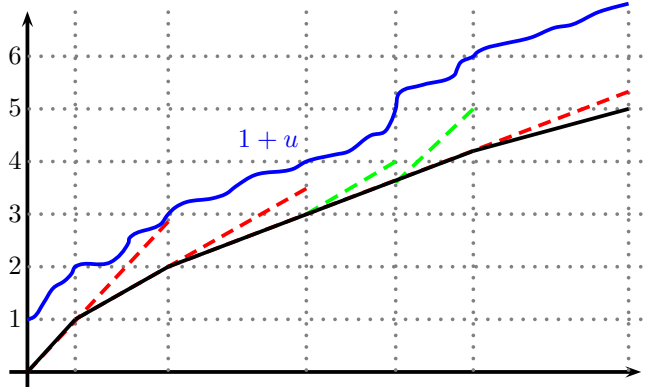
\begin{figure}
\begin{center}
\psset{xunit=1.2cm, yunit=1cm}
\psset{linewidth=0.5mm}
\begin{pspicture}(-.2,-.2)(6.8,4.8)
\psline{->}(0,-.2)(0,4.8)
\psline{->}(-.2,0)(6.8,0)
\psline[linestyle=dotted,linecolor=gray](0,0.7)(6.8,.7)\put(-0.25,0.6){1}
\psline[linestyle=dotted,linecolor=gray](0,1.4)(6.8,1.4)\put(-0.25,1.3){2}
\psline[linestyle=dotted,linecolor=gray](0,2.1)(6.8,2.1)\put(-0.25,2.0){3}
\psline[linestyle=dotted,linecolor=gray](0,2.8)(6.8,2.8)\put(-0.25,2.7){4}
\psline[linestyle=dotted,linecolor=gray](0,3.5)(6.8,3.5)\put(-0.25,3.4){5}
\psline[linestyle=dotted,linecolor=gray](0,4.2)(6.8,4.2)\put(-0.25,4.1){6}
\psline[linestyle=dotted,linecolor=gray](0.53,0)(0.53,4.8)
\psline[linestyle=dotted,linecolor=gray](1.56,0)(1.56,4.8)
\psline[linestyle=dotted,linecolor=gray](3.09,0)(3.09,4.8)
\psline[linestyle=dotted,linecolor=gray](4.08,0)(4.08,4.8)
\psline[linestyle=dotted,linecolor=gray](4.94,0)(4.94,4.8)
\psline[linestyle=dotted,linecolor=gray](6.66,0)(6.66,4.8)
\psline[linestyle=dashed,linecolor=red](0,0)(1.56,2.00)
\psline[linestyle=dashed,linecolor=red](0.53,0.7)(3.09,2.44)
\psline[linestyle=dashed,linecolor=red](1.56,1.4)(6.66,3.73)
\psline[linestyle=dashed,linecolor=green](3.09,2.1)(4.08,2.8)
\psline[linestyle=dashed,linecolor=green](4.08,2.52)(4.94,3.5)
\psline(0,0)(0.53,0.7)
\psline(0.53,0.7)(1.56,1.4)
\psline(1.56,1.4)(4.94,2.94)
\psline(4.94,2.94)(6.66,3.5)
\put(2.8,3){\blue$1+u$}
\listplot[showpoints=false,linecolor=blue,plotstyle=curve]{%
0.00 0.70 0.07 0.73 0.18 0.92 0.29 1.12 0.41 1.21 0.49 1.31 0.53 1.40
0.92 1.46 1.12 1.72 1.13 1.78 1.48 1.97 1.54 2.07 1.56 2.10 1.64 2.19
1.76 2.26 2.05 2.31 2.20 2.37 2.31 2.47 2.52 2.62 2.94 2.70 3.09 2.80
3.40 2.89 3.55 2.93 3.81 3.15 3.95 3.19 4.02 3.31 4.08 3.50 4.12 3.71
4.38 3.82 4.73 3.95 4.76 4.03 4.80 4.12 4.94 4.20 4.98 4.25 5.08 4.31
5.25 4.36 5.58 4.45 5.84 4.58 6.03 4.62 6.33 4.79 6.66 4.90}
\end{pspicture}
\caption{An example of the construction}
\end{center}
\end{figure}
All properties of $\phi$ are apparent from the picture, we only show that
$\phi(y_n)\uparrow\infty$. Let $A=\{n:\phi(y_n)=n\}$. If $A$ is infinite,
we are done, otherwise, let $N$ be the largest value in $A$, then
for $x\ge y_N$,
$$
\phi(x)=\phi(y_{N-1})+\frac{N-\phi(y_{N-1})}{y_N-y_{N-1}}(x-y_{N-1})
$$
and $\phi(x)\uparrow\infty$, since $\phi(x_{N-1})\le N-1<N$.
\end{proof}
\subsection{A slight variation of Gronwall's lemma}
Here we give a detailed proof of the variation of Gronwall's lemma used
in Section \ref{ss:exIM}. The result is elementary and probably well known,
it is given here only for the sake of completeness. The main differences are
the following: we do not assume that the term $a(\cdot)$ is positive
and the inequality holds only for a.\ e.\ time, but then it holds starting from
arbitrary initial times.
\begin{proposition}\label{p:gronwall}
Let $a$, $b\in L^1(0,T)$, with $b\ge0$ and let $u:[0,T]\to\mathbf{R}$ be
a lower semi-continuous and positive function. Assume that there exists
a set $S\subset (0,T]$ (thus, not containing $0$) with null Lebesgue
measure, such that for all $s\not\in S$ and all $t\in[s,T]$,
$$
u(t)\le u(s)+\int_s^ta(r)u(r)\,dr+\int_s^tb(r)\,dr.
$$
Then
$$
u(T)\le u(0)\,\mathrm{e}^{\int_0^Ta(s)\,ds}+\int_0^Tb(s)\,\mathrm{e}^{\int_s^Ta(r)\,dr}\,ds.
$$
\end{proposition}
\begin{proof}
We only need to prove the proposition if $a(\cdot)$ is piecewise constant.
Indeed, if this claim is true and $a\in L^1(0,T)$, there are piecewise
constant functions $a_n$ such that $a_n\longrightarrow a$ and without
loss of generality we can assume that each $a_n$ is constant on a finite
number of intervals whose extreme points do not belong to $S$ (but
possibly for the last one). By the usual Gronwall's lemma we can deduce
that $u$ is bounded by some constant $M$. We then set $b_n(s)=b(s)+M|a(s)-a_n(s)|$,
and we apply the claim with $a_n$ and $b_n$. As $n\to\infty$, we recover
the original statement.

Assume then that $a=\sum_{k=0}^{n-1}\alpha_k\mathbf{1}_{J_k}$,
where the intervals $J_k=[t_k,t_{k+1})$, $0=t_0<t_1<\dots<t_n=T$
and $t_0$, $t_1$, \ldots $t_{n-1}\not\in S$.
If $\alpha_k\ge0$, since $t_k\in S$, we know by the usual Gronwall's lemma
and semi-continuity of $u$ that
$$
u(t_{k+1})\le u(t_k)\,\mathrm{e}^{\alpha_k(t_{k+1}-t_k)}+\int_{t_k}^{t_{k+1}}b(s)\,\mathrm{e}^{\alpha_k(t_{k+1}-s)}\,ds.
$$
If $\alpha_k<0$, we reverse time as it is done in the proof of Theorem 5 of Flandoli
\& Romito \cite{FlRo01} and we apply again Gronwall's lemma to get
$$
u(t_{k+1})\le u(t_k)\,\mathrm{e}^{\alpha_k(t_{k+1}-t_k)}+\int_{t_k}^{t_{k+1}}b(s)\,\mathrm{e}^{\alpha_k(t_{k+1}-s)}\,ds.
$$
It is then easy to prove by induction on $k\le n$ that
$$
u(t_k)\le u(0)\,\mathrm{e}^{\int_0^{t_k}a(s)\,ds}+\int_0^{t_k}b(s)\,\mathrm{e}^{\int_s^{t_k}a(r)\,dr}\,ds.
$$
and in particular $k=n$ is exactly what we aimed to prove.
\end{proof}

\end{document}